\RequirePackage{ifthen}
\newboolean{MPA}
\setboolean{MPA}{false}

\ifthenelse {\boolean{MPA}}
{
\documentclass{svjour3}
\smartqed
\usepackage[margin=1in]{geometry}
} {
\documentclass[11pt]{article}
\usepackage[margin=.9in]{geometry}
}

\usepackage{graphicx}
\usepackage{amssymb}
\usepackage{amsmath}
\usepackage{verbatim,booktabs}
\usepackage{subfigure,epsfig,url,psfrag}
\usepackage{float}
\usepackage{mathrsfs}

\usepackage[usenames,dvipsnames]{pstricks}
\usepackage{epsfig}

\usepackage{cite} % to allow for spaces in citations

\usepackage{color}
\usepackage{enumerate}

\usepackage{todonotes}

\usepackage{hyperref}
\usepackage[capitalize,noabbrev]{cleveref}
\crefname{problem}{Problem}{Problems}
\crefname{claim}{Claim}{Claims}

\newcommand{\pare}[1]{\left(#1\right)}

\newcommand{\conv}{\mathop{\rm conv}}
\newcommand{\ie}{i.e., }

\renewcommand{\S}{\mathcal S}

\newcommand{\MP}{\text{MP}}
\newcommand{\LP}{\text{LP}}

\newcommand{\Q}{\mathcal Q}

\def\R{{\mathbb R}}

\def\01{\ensuremath{0\mathord{-}1}}
\allowdisplaybreaks

\ifthenelse {\boolean{MPA}}
{
\newenvironment{prf}[1][]
{\begin{proof}}
{\qed \end{proof}}

\newenvironment{prfc}[1][]
{\begin{proof}[#1]}
{\qed \end{proof}}

\newenvironment{prfh}[1][]
{\begin{proof}}
{\end{proof}}

\newcounter{claim} %[section]
\renewenvironment{claim}[1][]
{\refstepcounter{claim} \begin{trivlist} \item[] {\bf Claim~\theclaim}\space#1 \itshape}
{\end{trivlist}}

\newenvironment{cpf}
{\begin{trivlist} \item[] {\em Proof of claim }}
{$\hfill\diamond$ \end{trivlist}}

\journalname{Mathematical Programming A}

} {
\usepackage{amsthm}
\newtheorem{theorem}{Theorem}
\newtheorem{definition}{Definition}
\newtheorem{lemma}{Lemma}
\newtheorem{corollary}{Corollary}
\newtheorem{proposition}{Proposition}
\newtheorem{example}{Example}

\newtheorem{claim}{Claim}
\newtheorem{remark}{Remark}

\newenvironment{prf}[1][]
{\begin{proof}}
{\end{proof}}

\newenvironment{cpf}
{\begin{trivlist} \item[] {\em Proof of claim. }}
{$\hfill\diamond$ \end{trivlist}}

}

\usepackage{thmtools}
\usepackage{thm-restate}

\title{A polynomial-size extended formulation for the multilinear polytope of beta-acyclic hypergraphs}

\ifthenelse {\boolean{MPA}}
{
\titlerunning{A polynomial-size extended formulation for the multilinear polytope of beta-acyclic hypergraphs}
\author{Alberto Del Pia \and Aida Khajavirad}

\institute{Alberto Del Pia \at
              Department of Industrial and Systems Engineering \& Wisconsin Institute for Discovery, \\
              University of Wisconsin-Madison.
              E-mail: {\tt delpia@wisc.edu}.
           \and
           Aida Khajavirad \at
              Department of Industrial and Systems Engineering,
             Lehigh University.
             E-mail: {\tt aida@lehigh.edu}.
}
}
{
\author{Alberto Del Pia
\thanks{Department of Industrial and Systems Engineering \& Wisconsin Institute for Discovery,
             University of Wisconsin-Madison.
             E-mail: {\tt delpia@wisc.edu}.
             }
\and
Aida Khajavirad
\thanks{Department of Industrial and Systems Engineering,
             Lehigh University.
             E-mail: {\tt aida@lehigh.edu}.
             }
}
}

\date{August 29, 2023}

\begin{document}

\maketitle

\begin{abstract}
We consider the multilinear polytope defined as the convex hull of the set of binary points $z$, satisfying a collection of equations of the form
$z_e = \prod_{v \in e} z_v$ for all $e \in E$. The complexity of the facial structure of the multilinear polytope is closely related to the acyclicity degree of the underlying hypergraph. We obtain a polynomial-size extended formulation for the multilinear polytope of $\beta$-acyclic hypergraphs, hence characterizing the acyclic hypergraphs for which such a formulation can be constructed.
%, hence, settling completely the extension complexity of the multilinear polytope of acyclic hypergraphs.
\ifthenelse {\boolean{MPA}}
{
% For MPA begin
% For MPA end
} {}
\end{abstract}

\ifthenelse {\boolean{MPA}}
{}{
% For OO begin
}

\section{Introduction}
Binary polynomial optimization, \ie the problem of finding a binary point maximizing a polynomial function, is a fundamental NP-hard problem in discrete optimization with a wide range of applications across science and engineering. To formally define this problem, we employ a hypergraph representation scheme  introduced in~\cite{dPKha17MOR}. A \emph{hypergraph} $G$ is a pair $(V,E)$, where $V$ is a finite set of nodes and $E$ is a set of subsets of $V$, called the edges of $G$. Throughout this paper we consider hypergraphs without loops or parallel edges, in which case $E$ is a set of subsets of $V$ of cardinality at least two. Moreover, the \emph{rank} of a hypergraph $G$ is the maximum cardinality of any edge in $E$. 
With any hypergraph $G= (V,E)$, we associate the following binary polynomial optimization problem:
\begin{align}
\label[problem]{prob BP}
\tag{BP}
\begin{split}
\max & \qquad \sum_{v\in V} {c_v z_v} + \sum_{e\in E} {c_e \prod_{v\in e} {z_v}} \\
{\rm s.t.} & \qquad z_v \in \{0,1\} \qquad \forall v \in V,
\end{split}
\end{align}
where without loss of generality we assume $c_e \neq 0$ for all $e \in E$.  
%With the objective of constructing convex relaxations for Problem~\eqref{prob BP}, 
Following a common practice in nonconvex optimization,  we then proceed with linearizing the objective function by introducing new variables for each product term to obtain an equivalent reformulation of \cref{prob BP} in a lifted space:
\begin{align}
\label[problem]{prob LBP}
\tag{LBP}
\begin{split}
\max & \qquad \sum_{v\in V} {c_v z_v} + \sum_{e\in E} {c_e z_e} \\
{\rm s.t.} & \qquad z_e = \prod_{v\in e} {z_v} \qquad \forall e \in E\\
& \qquad z_v \in \{0,1\} \qquad \forall v \in V.
\end{split}
\end{align}

\subsection{The multilinear polytope and hypergraph acyclicity} 

To solve \cref{prob LBP} efficiently using polyhedral techniques, it is essential to understand the facial structure of the polyhedral convex hull of its feasible region. 
To this end, in the same vein as~\cite{dPKha17MOR}, we define the~\emph{multilinear set} as
\begin{equation*} \label{eq: SG}
\S_G= \Big\{ z \in \{0,1\}^{V+E} : z_e = \prod_{v \in e} {z_{v}}, \; \forall e \in E \Big\},
\end{equation*}
and we refer to its convex hull as the~\emph{multilinear polytope} and denote it by $\MP_G$. A simple polyhedral relaxation of $\S_G$ can be obtained by replacing each term $z_e = \prod_{v \in e} {z_{v}}$ by its convex hull over the unit hypercube: 
\begin{equation*}
%\label{stdre}
\MP^{\LP}_G =\Big\{z: z_v \leq 1, \forall v \in V,\; z_e  \geq 0, \; z_e \geq \sum_{v\in e}{z_v}-|e|+1, \forall e \in E, \;
z_e  \leq z_v, \forall e \in E, \forall v \in e\Big\}.
\end{equation*}
The above relaxation is often referred to as the \emph{standard linearization} and has been used extensively in the literature~\cite{yc93}. 
In the special case with $r = 2$; \ie when all product terms in $\S_G$ are products of two variables, the multilinear polytope coincides with the well-known Boolean quadric polytope ${\rm BQP}_G$~\cite{Pad89}. Padberg~\cite{Pad89} proves that ${\rm BQP}_G$ coincides with its standard linearization if and only if the graph $G$ is acyclic. Hence it is natural to ask whether the multilinear polytope of acyclic hypergraphs has a simple structure as well.
Unlike graphs, the notions of cycles and acyclicity in hypergraphs are not unique. The most well-known types of acyclic hypergraphs, in increasing order of generality, are Berge-acyclic, $\gamma$-acyclic, $\beta$-acyclic, and $\alpha$-acyclic hypergraphs~\cite{fagin83,BeFaMaYa83,Dur12,bra14}. In the following, we present a brief review of the literature on the mutlilinear polytope of acyclic hypergraphs.  

%In~\cite{dPKha18SIOPT,dPKha21MOR}, the authors relate the complexity of the multilinear polytope to the \emph{acyclicity degree} of its hypergraph.
In~\cite{dPKha18SIOPT,BucCraRod16}, the authors prove that $\MP_G = \MP^{\LP}_G$ if and only if the hypergraph $G$ is Berge-acyclic.  
In~\cite{dPKha18SIOPT}, the authors introduce \emph{flower inequalities}, a class of facet-defining inequalities for the multilinear polytope, and show that the polytope obtained by adding all such inequalities to $\MP^{\LP}_G$ 
coincides with $\MP_G$ if and only if the hypergraph $G$ is $\gamma$-acyclic.
While the multilinear polytope of $\gamma$-acyclic hypergraphs may contain exponentially many facets, a polynomial-size extended formulation of $\MP_G$ is implicit in~\cite{dPKha18SIOPT}.
\footnote{By~\emph{polynomial-size extended formulation,} we mean that the size of the system of linear inequalities, as defined in \cite{SchBookIP}, is polynomial in the number of nodes and edges of $G$, which is a stronger notion than asking for a polynomial number of variables and inequalities.}
%the authors present a strongly polynomial-time algorithm to solve the separation problem.
Subsequently, in~\cite{dPKha21MOR}, the authors  introduce~\emph{running intersection inequalities}, a class of facet-defining inequalities for the multilinear polytope that serve as a generalization of flower inequalities. The authors prove that
for kite-free $\beta$-acyclic hypergraphs, a class that lies
between $\gamma$-acyclic and $\beta$-acyclic hypergraphs, the polytope obtained by adding all running intersection inequalities to $\MP^{\LP}_G$ coincides with $\MP_G$, and it admits a polynomial-size extended formulation. 
At the other end of the spectrum, in~\cite{delGre22,dPDiG23ALG}, the authors prove that \cref{prob BP} is strongly NP-hard over $\alpha$-acyclic hypergraphs.
This result implies that, unless P = NP, one cannot construct a polynomial-size extended formulation for the multilinear polytope of $\alpha$-acyclic hypergraphs. 
See~\cite{CraRod16,BieMun18,HojPfeWal19,dPKhaSah20MPC,chenSanOkt20,XuAdaAks20,dPDiG21IJO,Aida22,kim22,dPWal22IPCO} for further results regarding polyhedral relaxations of multilinear sets.

Hence, to this date, there remains one class of acyclic hypergraphs for which we do not know whether it is possible to obtain a polynomial-size extended formulation: the class of $\beta$-acyclic hypergraphs.
% Hence, to this date, there remains one last open question regarding the extension complexity 
% %facial structure 
% of the multilinear polytope of acyclic hypergraphs: does the mutilinear polytope of $\beta$-acyclic hypergraphs admit a  polynomial-size extended formulation?
%Recall that the \emph{extension complexity} of a polytope $\P$ is defined as the smallest number of inequalities necessary to describe a higher dimensional polytope $\Q$ that can be projected on $\P$. 
%
In~\cite{delGre22,dPDiG23ALG}, the authors present a strongly polynomial time algorithm to solve \cref{prob BP} over $\beta$-acyclic hypergraphs. 
While this result settles the algorithmic complexity of \cref{prob BP} over acyclic hypergraphs, it does not address the complexity of the extended formulation.
%settle the extension complexity of the multilinear polytope of $\beta$-acyclic hypergraphs.
%In fact, in the celebrated paper \cite{rothvoss17}, Rothvoss showed that 
Indeed, it is well-known that there exist polytopes over which one can optimize any linear function in strongly polynomial time, yet they do not admit any polynomial-size extended formulation (see, e.g., \cite{rothvoss17}).

%See for example a celebrated result of Rothvoss~\cite{rothvoss17}, which implies that there exists no polynomial-size extended formulation for the matching polytope, even though any linear function can be maximized over the matching polytope in strongly polynomial time using Edmond's Blossom algorithm. 

We should remark that it is possible and in fact highly plausible that there exists a family of hypergraphs between $\alpha$-acyclic and $\beta$-acyclic hypergraphs for which one can obtain a polynomial-size extended formulation of the multilinear polytope. However, our focus in this paper is to characterize the known classes of acyclic hypergraphs for which it is possible to construct a polynomial-size extended formulation.

\subsection{Our contribution}

In this paper, we present a polynomial-size extended formulation for the multilinear polytope of $\beta$-acyclic hypergraphs. Recall that a \emph{$\beta$-cycle} of length $t$, for some $t \geq 3$, is a sequence $v_1, e_1, v_2, e_2, \dots, v_t, e_t, v_1$ such that $v_1$, $v_2$, $\dots$, $v_t$ are distinct nodes, $e_1$, $e_2$, $\dots$, $e_t$ are distinct edges, and $v_i$ belongs to $e_{i-1}, e_i$ and no other $e_j$, for all $i = 1,\dots,t$, where we define $e_0 := e_t$.
A hypergraph is called~\emph{$\beta$-acyclic} if it does not contain any $\beta$-cycle.  The following statement summarizes our main result regarding the existence of a polynomial-size extended formulation for the multilinear polytope of $\beta$-acyclic hypergraphs:

\begin{restatable}{theorem}{First}
%\begin{theorem} 
\label{th main}
Let $G = (V,E)$ be a $\beta$-acyclic hypergraph of rank $r$.
Then there exists an polynomial-size extended formulation of $\MP_G$ comprising of at most 
$(3r-4)|V|+4|E|$ inequalities, with at most $(r-2) |V|$ extended variables.
The system is explicitly given in~\cref{th2}.
%\end{theorem}
\end{restatable}

%The proof of Theorem~\ref{th main} is given in Section~xxx.

%More precisely, given a $\beta$-acyclic hypergraph $G =(V, E)$ of rank $r$, we present an extended formulation for $\MP_G$ comprising of at most $(3r-4)|V|+4|E|$ inequalities, and with at most $(r-2)|V|$ extra variables. 

It is important to note that the standard linearization of a rank $r$ hypergraph $G=(V, E)$ consists of $2|V|+(r+2) |E|$ inequalities. It is well-understood that $\MP^{\rm LP}_G$ often leads to very weak relaxations of $\MP_G$ for $\beta$-acyclic hypergraphs. Theorem~\ref{th main} implies that while the proposed extended formulation for $\MP_G$ contains $(r-2)|V|$ additional variables, it has fewer inequalities than the standard linearization for $\beta$-acyclic hypergraphs with $|E| \geq 3 |V|$. We should also remark that the inequalities defining our proposed extended formulation are very sparse; that is, they contain at most four variables with non-zero coefficients; a feature that is highly beneficial from a computational perspective.

%\aida{Once you decide how to make this change lets implement it together so that i tell you where to move stuff in the main text?}
Our construction relies on the key concept of nest points of hypergraphs. A node $v \in V$ is a \emph{nest point} of $G$ if the set of the edges of $G$ containing $v$ is totally ordered. In other words, the edges in $E$ containing $v$ can be ordered so that $e_1 \subset e_2 \subset \cdots \subset e_k$.
%if for every two edges $e,f \in E$ containing $v$, either $e \subseteq f$ or $f \subseteq e$.
It is simple to see that nest points can be found in polynomial time. We define the hypergraph obtained from $G = (V,E)$ by \emph{removing} a node $v \in V$ as $G - v := (V',E')$, where $V' := V \setminus \{ v \}$ and $E' := \{ e \setminus \{v\} : e \in E, \ |e \setminus \{v\}| \ge 2 \}$. 
A \emph{nest point sequence} of length $s$ for some  $s\le |V|$ of $G$ is an ordering $v_1, \dots, v_s$ of $s$ distinct nodes of $G$, such that $v_1$ is a nest point of $G$, $v_2$ is a nest point of $G - v_1$, and so on, until $v_s$ is a nest point of $G - v_1 - \dots - v_{s-1}$.
We can write this condition compactly as $v_i$ is a nest point of $G - v_1 - \dots - v_{i-1}$, for $i=1,\dots,s$, 
where we make the slight abuse of notation $G - v_1 - \dots - v_0 = G$.  We then use the following characterization of $\beta$-acyclic hypergraphs, in terms of nest points:
%which relates $\beta$-acyclicity to nest points:
%We are now ready to state this characterization of $\beta$-acyclic hypergraphs.
\begin{theorem}[\cite{Dur12}]\label{th beta iff}
A hypergraph $G$ is $\beta$-acyclic if and only if after removing recursively a nest point, until one is found, we obtain the empty hypergraph $(\emptyset,\emptyset)$.
\end{theorem}

From \cref{th beta iff} it follows that a hypergraph is $\beta$-acyclic if and only if it has a nest point sequence of length $|V|$. 
In fact, our approach to prove \cref{th main} can be used to obtain extended formulations for the multilinear polytope of  more general hypergraphs containing $\beta$-cycles; namely, hypergraphs containing a nest point sequence of length $s$ for some $1 \leq s \leq |V|$:

%\begin{theorem} 
\begin{restatable}{theorem}{Second}
\label{int2}
Let $G = (V,E)$ be a hypergraph of rank $r$, and let $v_1,\dots,v_s$ be a nest point sequence of $G$.
Then an extended formulation of $\MP_G$ is given by a description of $\MP_{G - v_1 - \cdots - v_s}$, together with a system of at most $|V|+2|E| + 4 rs$ linear inequalities, including at most $(r-2) s$ extended variables.
The system is characterized in \cref{th general s nodes}.
\end{restatable}

%The proof of Theorem~\ref{int2} is given in Section~xxx.

To prove Theorems~\ref{th main} and~\ref{int2}, we present, in \cref{th decomp}, a new sufficient condition for decomposability of multilinear sets that is of independent interest.

A natural question is whether it is possible to characterize the multilinear polytope of $\beta$-acyclic hypergraphs in the original space of variables.
We argue that for a $\beta$-acyclic hypergraph $G$, an explicit description of $\MP_G$ in the original space does not have desirable numerical properties, as this polytope may contain 
\emph{very dense} facet-defining inequalities. To demonstrate this property, in \cref{prop dense} we present a family  of $\beta$-acyclic hypergraphs $G=(V,E)$
whose multilinear polytope consists of facet-defining inequalities 
containing $|E|$ variables with non-zero coefficients.
It is well-understood that the addition of such dense inequalities as cutting planes to an LP relaxation in a branch-and-cut solver often leads to increased CPU times.
Finally, as a byproduct of our convex hull characterizations, we present a new class of sparse valid inequalities for the multilinear polytope in the original space, which serve as a generalization of running intersection inequalities \cite{dPKha21MOR}. These inequalities can be incorporated in branch-and-cut based global solvers to improve the quality of existing relaxations for nonconvex problems whose factorable reformulations contain multilinear sets~\cite{IdaNick18,dPKhaSah20MPC}.

\paragraph{Outline.} The remainder of this paper is organized as follows. In Section~\ref{sec: decomp}, we present a sufficient condition for decomposability of multilinear sets that enables us to decompose multilinear sets of hypergraphs with nest points to simpler multilinear sets (see Theorem~\ref{th decomp}).
In Section~\ref{sec: pointed}, we consider a special type of hypergraphs obtained as a result of decomposing hypergraphs with nest points, and characterize its multilinear polytope using a direct approach (see Theorem~\ref{th pointed MPG}). In Section~\ref{sec: extended}, by combining the results of Sections~\ref{sec: decomp} and~\ref{sec: pointed}, we describe the mulilinear poyltope of hypergraphs with nest points in terms of multilinear polytopes of simpler hypergraphs (see Theorems~\ref{int2} and ~\ref{th general s nodes}). Subsequently, we obtain a polynomial-size extended formulation for the multilinear polytope of $\beta$-acyclic hypergraphs (see Theorems~\ref{th main} and~\ref{th2}). In Section~\ref{sec:original}, we elaborate on the complexity of the multilinear polytope of $\beta$-acyclic hypergraphs in the original space. We conclude by presenting a new class of sparse valid inequalities for the multilinear polytope of general hypergraphs.

\section{Decomposability of multilinear sets}
\label{sec: decomp}

In this section, we present a new sufficient condition for decomposability of multilinear sets that we will use to obtain our extended formulations for the multilinear polytope of hypergraphs with nest points.

Consider hypergraphs $G_1 = (V_1,E_1)$ and $G_2=(V_2,E_2)$ such that $V_1 \cap V_2 \neq \emptyset$.
We denote by $G_1 \cap G_2$ the hypergraph $(V_1 \cap V_2, E_1 \cap E_2)$ and 
by $G_1 \cup G_2$ the hypergraph $(V_1 \cup V_2, E_1 \cup E_2)$.
Let $G := G_1 \cup G_2$.
We say that the set $\S_G$ is \emph{decomposable into the sets $\S_{G_1}$ and $\S_{G_2}$} if
\begin{align*}
\conv\S_{G} = \conv \bar \S_{G_1} \cap \conv \bar \S_{G_2},
\end{align*}
where $\bar \S_{G_1}$ (resp.~$\bar \S_{G_2}$) is the set of all points in the space of $\S_G$ whose projection in the space defined by $G_1$ (resp.~$G_2$) is $\S_{G_1}$ (resp.~$\S_{G_2}$). 

Other known decomposition results for multilinear sets are Theorem~1 in~\cite{dPKha18MPA}, % and 4
Theorem~5 in~\cite{dPKha18SIOPT},
Theorem~1 in~\cite{dPKha21MOR},
and Theorem~4 in~\cite{dPDiG21IJO}.
In all prior decomposition results, the hypergraphs $G_1$ and $G_2$ are assumed to be section hypergraphs of $G$. 
Recall that $G_1$ is a \emph{section hypergraph} of $G=(V,E)$ if $G_1 = (V_1,E_1)$, where $V_1 \subset V$ and $E_1 = \{e \in E: e \subseteq V_1\}$.
This means that $G_1$ and $G_2$ inherit all edges of $G$ contained in their respective node sets. 
On the contrary, in our new decomposition result, $G_1$ is generally not a section hypergraph of $G$, and this key difference allows $G_1$ to have a very simple structure that will be exploited in \cref{sec: pointed}.

In the remainder of the paper, for notational simplicity, we sometimes write a node variable $z_v$ as $z_{\{v\}}$.
This can happen, for example, when we have an edge $e$ of cardinality two, $e = \{u,v\}$, and we write the variable corresponding to $u$ as $z_{e \setminus \{v\}} = z_{\{u\}}$.
%
%\medskip
We now present our decomposition result.

\begin{theorem}
\label{th decomp}
Let $G = (V,E)$ be a hypergraph, let $v$ be a nest point of $G$, let $e_1 \subset e_2 \subset \cdots \subset e_k$ be the edges of $G$ containing $v$, and let $E_v := \{e_1,\dots,e_k\}$.
For each $i \in [k]:= \{1,\cdots,k\}$, let $p_i := e_i \setminus \{v\}$ and define $P_v := \{p \in \{p_1,\dots,p_k\} : |p| \ge 2\}$.
Assume that $P_v \subseteq E$.
Let $G_1 := (e_k, E_v \cup P_v)$ and let $G_2 := G - v$.
Then the set $\S_G$ is decomposable into $\S_{G_1}$ and $\S_{G_2}$.
\end{theorem}

\begin{prf}
We assume $k \ge 1$, as otherwise the result is obvious.
%Denote by $e_1,\dots,e_k$, for $k \ge 1$, the edges in $E_v$.
%Since $v$ is a nest point, without loss of generality we assume $e_1 \subset e_2 \subset \cdots \subset e_k$.
%Let $\ell$ be the smallest index $i \in [k]$ with $|e_i \setminus \{v\}| \ge 2$.
%For $i \in [k]$, we define $p_i := e_i \setminus \{v\}$.
%Note that $P_v = \{p_1,\dots,p_k\}$, $p_1 \subset p_2 \subset \cdots \subset p_k$, and $p_1$ may have cardinality one.
%We also assume that $P_v \neq \emptyset$, otherwise node $v$ is contained in only one edge of $G$ and such edge has cardinality two, thus the result follows from theorem~1 in \cite{dPKha18MPA}. 
%
We now explain how we write, in the rest of the proof, a vector $z$ in the space defined by $G$ by partitioning its components in a number of subvectors.
The vector $z_\cap$ contains the components of $z$ corresponding to nodes and edges in $G_1 \cap G_2$, i.e., nodes in $e_k \setminus \{v\}$ and edges in $P_v$.
The vector $z_1$ contains the components of $z$ corresponding to nodes and edges in $G_1$ but not in $G_2$, i.e., node $v$ and edges in $E_v$.
Finally, the vector $z_2$ contains the components of $z$ corresponding to nodes and edges in $G_2$ but not in $G_1$.
Using these definitions, we can now write, up to reordering variables, $z=(z_1,z_\cap,z_2)$.
Similarly, we can write a vector $z$ in the space defined by $G_1$ as $(z_1,z_\cap)$, and a vector $z$ in the space defined by $G_2$ as $z=(z_\cap,z_2)$.

We now proceed with the proof of the theorem.
To this end, we show the two inclusions $\conv\S_{G} \subseteq \conv \bar \S_{G_1} \cap \conv \bar \S_{G_2}$ and $\conv\S_{G} \supseteq \conv \bar \S_{G_1} \cap \conv \bar \S_{G_2}$.
The first inclusion clearly holds, since $\S_G \subseteq \bar \S_{G_1} \cap \bar \S_{G_2}$. 
Thus, it suffices to show the inclusion $\conv\S_{G} \supseteq \conv \bar \S_{G_1} \cap \conv \bar \S_{G_2}$.
Let $\tilde z \in \conv \bar \S_{G_1} \cap \conv \bar \S_{G_2}$.
We will show that $\tilde z \in \conv\S_G$.

To prove $\tilde z \in \conv\S_G$, we will write $\tilde z$ explicitly as a convex combinations of vectors in $\S_G$.
In the next claim, we show how a vector in $\S_{G_1}$ and a vector in $\S_{G_2}$ can be combined to obtain a vector in $\S_G$.

\begin{claim}
\label{claim combination}
Let $(z_1,z_\cap) \in \S_{G_1}$ and $(z'_\cap,z'_2) \in \S_{G_2}$ such that 
$z_{p_i} = z'_{p_i}$ for every $i \in [k]$.
Then, $(z_1,z'_\cap,z'_2) \in \S_G$.
\end{claim}

\begin{cpf}
It suffices to show that $(z_1,z'_\cap) \in \S_{G_1}$.
The edges of $G_1$ whose components are in $z'_\cap$ are the edges in $P_v$.
Each edge in $P_v$ contains only nodes with components in $z'_\cap$, thus we have $z'_{p_i} = \prod_{u \in p_i} z'_u$, for each $p_i \in P_v$.
The edges of $G_1$ whose components are in $z_1$ are the edges in $E_v$, thus we only need to show $z_{e_i} = z_v \prod_{u \in p_i} z'_u$, for each $i \in [k]$.
This equality holds since
$$
z_{e_i} = z_v z_{p_i} = z_v z'_{p_i} = z_v \prod_{u \in p_i} z'_u.
$$
\end{cpf}

In the remainder of the proof, we show how to write explicitly $\tilde z$ as a convex combination of the vectors in $\S_G$ obtained in \cref{claim combination}.

By assumption, the vector $(\tilde z_1,\tilde z_\cap)$ is in $\conv \S_{G_1}$. Thus, it can be written as a convex combination of points in $\S_{G_1}$;
\ie there exists $\mu \ge 0$ such that
\begin{align}
%\label{eq mu sum}
\nonumber
\sum_{(z_1 ,z_\cap) \in \S_{G_1}} \mu_{(z_1,z_\cap)} & = 1 \\
\label{eq mu conv}
\sum_{(z_1 ,z_\cap) \in \S_{G_1}} \mu_{(z_1,z_\cap)} (z_1,z_\cap) & = (\tilde z_1,\tilde z_\cap).
\end{align}
Similarly, the vector $(\tilde z_\cap,\tilde z_2)$ is in $\conv \S_{G_2}$ and it can be written as a convex combination of points in $\S_{G_2}$;
\ie there exists $\nu \ge 0$ such that 
\begin{align}
%\label{eq nu sum}
\nonumber
\sum_{(z_\cap, z_2) \in \S_{G_2}} \nu_{(z_\cap,z_2)} & = 1 \\
\label{eq nu conv}
\sum_{(z_\cap, z_2) \in \S_{G_2}} \nu_{(z_\cap,z_2)} (z_\cap,z_2) & = (\tilde z_\cap,\tilde z_2).
\end{align}

For ease of notation, we define, for $i \in [k+1]$,
\begin{align*}
m(i) := 
\begin{cases}
1 - \tilde z_{p_1}  & \qquad \text{if $i = 1$} \\
\tilde z_{p_{i-1}} - \tilde z_{p_{i}} & \qquad \text{if $i \in \{2,\dots,k\}$} \\
\tilde z_{p_k}  & \qquad \text{if $i = k+1$}.
\end{cases}
\end{align*}
In the remainder of the proof, given binary $z_{p_1}, \dots, z_{p_k}$, we will consider the number $\min\{j \in [k+1] : z_{p_j} = 0\} \in \{1,\dots,k+1\}$, with the understanding that this number equals $k+1$ when $z_{p_1} = \cdots = z_{p_k} = 1$.
For every $(z_1,z_\cap) \in \S_{G_1}$ and $(z'_\cap,z'_2) \in \S_{G_2}$ such that
$z_{p_i} = z'_{p_i}$ for every $i \in [k]$, 
we define 
\begin{align*}
\lambda_{(z_1,z'_\cap,z'_2)} :=
\frac{\mu_{(z_1,z_\cap)} \nu_{(z'_\cap,z'_2)}}
{m(i)},
\end{align*}
where $i := \min\{j \in [k+1] : z_{p_j} = 0\} = \min\{j \in [k+1] : z'_{p_j} = 0\}$.
In the next claims we show that the vector $\lambda$ that we just defined serves as the vector of multipliers to write $\tilde z$ as a convex combination of the vectors in $\S_G$ obtained in \cref{claim combination}.
We start with a technical claim.
\begin{claim}
\label{claim sum mu nu}
For $i \in [k+1]$, we have
\begin{align*}
\sum_{\substack{(z_1 ,z_\cap) \in \S_{G_1} \\ \min\{j \in [k+1] : z_{p_j} = 0\} = i}} \mu_{(z_1,z_\cap)} 
=
\sum_{\substack{(z_\cap, z_2) \in \S_{G_2} \\ \min\{j \in [k+1] : z_{p_j} = 0\} = i}} \nu_{(z_\cap, z_2)} 
=
m(i)
\end{align*}
\end{claim}

\begin{cpf}
By considering the component of \eqref{eq mu conv} corresponding to $p_i$, for $i \in [k]$, we obtain
\begin{align*}
\sum_{\substack{(z_1 ,z_\cap) \in \S_{G_1} \\ z_{p_i}=1}} \mu_{(z_1,z_\cap)} 
%& = \sum_{\substack{(z_\cap, z_2) \in \S_{G_2} \\ z_{p_i}=1}} \nu_{(z_\cap,z_2)} 
= \tilde z_{p_i}.
%\sum_{\substack{(z_1 ,z_\cap) \in \S_{G_1} \\ z_{p_i}=0}} \mu_{(z_1,z_\cap)} 
%& = \sum_{\substack{(z_\cap, z_2) \in \S_{G_2} \\ z_{p_i}=0}} \nu_{(z_\cap,z_2)}
%= 1 - \tilde z_{p_k}.
\end{align*}

%We now get back to the statement of the claim, and consider the case $i=k+1$.
%We have
We first consider the case $i=k+1$.
We have
\begin{align*}
\sum_{\substack{(z_1 ,z_\cap) \in \S_{G_1} \\ \min\{j \in [k+1] : z_{p_j} = 0\} = k+1}} \mu_{(z_1,z_\cap)} 
=
\sum_{\substack{(z_1 ,z_\cap) \in \S_{G_1} \\ z_{p_k}=1}} \mu_{(z_1,z_\cap)} 
= \tilde z_{p_k}.
\end{align*}

Next, we consider the case $i=1$.
We have
\begin{align*}
\sum_{\substack{(z_1 ,z_\cap) \in \S_{G_1} \\ \min\{j \in [k+1] : z_{p_j} = 0\} = 1}} \mu_{(z_1,z_\cap)} 
=
\sum_{\substack{(z_1 ,z_\cap) \in \S_{G_1} \\ z_{p_1}=0}} \mu_{(z_1,z_\cap)} 
=
1 - \sum_{\substack{(z_1 ,z_\cap) \in \S_{G_1} \\ z_{p_1}=1}} \mu_{(z_1,z_\cap)} 
=
1 - \tilde z_{p_1}.
\end{align*}

% The case $i=k$:
% \begin{align*}
% \sum_{\substack{(z_1 ,z_\cap) \in \S_{G_1} \\ \min\{j \in [k+1] : z_{p_j} = 0\} = k}} \mu_{(z_1,z_\cap)} 
% =
% \sum_{\substack{(z_1 ,z_\cap) \in \S_{G_1} \\ z_{p_{k-1}}=1 \\ z_{p_{k}}=0}} \mu_{(z_1,z_\cap)} 
% =
% \sum_{\substack{(z_1 ,z_\cap) \in \S_{G_1} \\ z_{p_{k-1}}=1}} \mu_{(z_1,z_\cap)} 
% -
% \sum_{\substack{(z_1 ,z_\cap) \in \S_{G_1} \\ \min\{j \in [k+1] : z_{p_j} = 0\} = k+1}}
% %z_{p_{k}}=1}} 
% \mu_{(z_1,z_\cap)}  
% = \tilde z_{p_{k-1}} - \tilde z_{p_k}.
% \end{align*}

% The case $i=k-1$:
% \begin{align*}
% \sum_{\substack{(z_1 ,z_\cap) \in \S_{G_1} \\ \min\{j \in [k+1] : z_{p_j} = 0\} = k-1}} \mu_{(z_1,z_\cap)} 
% & =
% \sum_{\substack{(z_1 ,z_\cap) \in \S_{G_1} \\ z_{p_{k-2}}=1 \\ z_{p_{k-1}}=0}} \mu_{(z_1,z_\cap)} \\
% & =
% \sum_{\substack{(z_1 ,z_\cap) \in \S_{G_1} \\ z_{p_{k-2}}=1}} \mu_{(z_1,z_\cap)} 
% -
% \sum_{\substack{(z_1 ,z_\cap) \in \S_{G_1} \\ \min\{j \in [k+1] : z_{p_j} = 0\} = k}} \mu_{(z_1,z_\cap)}  
% -
% \sum_{\substack{(z_1 ,z_\cap) \in \S_{G_1} \\ \min\{j \in [k+1] : z_{p_j} = 0\} = k+1}} \mu_{(z_1,z_\cap)} \\
% &
% = \tilde z_{p_{k-2}} - (\tilde z_{p_{k-1}} - \tilde z_{p_k}) - \tilde z_{p_k}
% = \tilde z_{p_{k-2}} - \tilde z_{p_{k-1}}.
% \end{align*}

Next, we consider the case $i \in [k]$.
We have
\begin{align*}
\sum_{\substack{(z_1 ,z_\cap) \in \S_{G_1} \\ \min\{j \in [k+1] : z_{p_j} = 0\} = i}} \mu_{(z_1,z_\cap)} 
& =
\sum_{\substack{(z_1 ,z_\cap) \in \S_{G_1} \\ z_{p_{i-1}}=1 \\ z_{p_{i}}=0}} \mu_{(z_1,z_\cap)} \\
& =
\sum_{\substack{(z_1 ,z_\cap) \in \S_{G_1} \\ z_{p_{i-1}}=1}} \mu_{(z_1,z_\cap)} 
-
\sum_{\substack{(z_1 ,z_\cap) \in \S_{G_1} \\ z_{p_{i}}=1}} \mu_{(z_1,z_\cap)} \\
&
= \tilde z_{p_{i-1}} - \tilde z_{p_{i}}.
\end{align*}

The statement for $\nu$ follows symmetrically, starting with \eqref{eq nu conv} rather than \eqref{eq mu conv}.
\end{cpf}

In the next claim, we show that the multipliers $\lambda$ are nonnegative and sum to one.

\begin{claim}
\label{claim sum one}
We have $\lambda \ge 0$ and 
\begin{align*}
\sum_{\substack{(z_1 ,z_\cap) \in \S_{G_1} \\ (z'_\cap, z'_2) \in \S_{G_2} \\ z_{p_i} = z'_{p_i} \forall i \in [k]}} 
\lambda_{(z_1,z'_\cap,z'_2)} =1
\end{align*}
\end{claim}

\begin{cpf}
% It is simple to observe that $p_1 \subset p_2 \subset \cdots \subset p_k$ implies 
% $\tilde z_{p_1} \le \tilde z_{p_2} \le \cdots \le z_{p_k}$, thus $m(i) \ge 0$
It follows from \cref{claim sum mu nu} that $m(i) \ge 0$ for all $i \in [k+1]$.
Thus, using the definition of $\lambda$, we obtain $\lambda \ge 0$.
%
%Using the fact that $p_1 \subset p_2 \subset \cdots \subset p_k$, 
Using \cref{claim sum mu nu}, we obtain
\begin{align*}
\sum_{\substack{(z_1 ,z_\cap) \in \S_{G_1} \\ (z'_\cap, z'_2) \in \S_{G_2} \\ z_{p_i} = z'_{p_i} \forall i \in [k]}} 
\lambda_{(z_1,z'_\cap,z'_2)} 
& = 
\sum_{i \in [k+1]} \ 
\sum_{\substack{(z_1 ,z_\cap) \in \S_{G_1} \\ \min\{j \in [k+1] : z_{p_j} = 0\} = i}} \
\sum_{\substack{(z'_\cap, z'_2) \in \S_{G_2} \\ \min\{j \in [k+1] : z'_{p_j} = 0\} = i}} \
\lambda_{(z_1,z'_\cap,z'_2)} \\
& = 
\sum_{i \in [k+1]} \ 
\sum_{\substack{(z_1 ,z_\cap) \in \S_{G_1} \\ \min\{j \in [k+1] : z_{p_j} = 0\} = i}} \
\sum_{\substack{(z'_\cap, z'_2) \in \S_{G_2} \\ \min\{j \in [k+1] : z'_{p_j} = 0\} = i}} \
\frac{\mu_{(z_1,z_\cap)} \nu_{(z'_\cap,z'_2)}}{m(i)} \\
& = 
\sum_{i \in [k+1]} 
\frac{1}{m(i)}
\pare{
\sum_{\substack{(z_1 ,z_\cap) \in \S_{G_1} \\ \min\{j \in [k+1] : z_{p_j} = 0\} = i}} 
\mu_{(z_1,z_\cap)}}
\pare{
\sum_{\substack{(z'_\cap, z'_2) \in \S_{G_2} \\ \min\{j \in [k+1] : z'_{p_j} = 0\} = i}} 
\nu_{(z'_\cap,z'_2)} } \\
& = 
\sum_{i \in [k+1]} 
\frac{(m(i))^2}{m(i)} 
= 
\sum_{i \in [k+1]} 
m(i) 
= 1.
\end{align*}
\end{cpf}

Our last claim, which concludes the proof of the theorem, shows that the multipliers $\lambda$ yield $\tilde z$ as a convex combination of the vectors in $\S_G$ obtained in \cref{claim combination}.

\begin{claim}
\label{claim conv comb}
We have 
\begin{align}
\label{eq to verify}
(\tilde z_1,\tilde z_\cap,\tilde z_2) = 
\sum_{\substack{(z_1 ,z_\cap) \in \S_{G_1} \\ (z'_\cap, z'_2) \in \S_{G_2} \\ z_{p_i} = z'_{p_i} \forall i \in [k]}} 
\lambda_{(z_1,z'_\cap,z'_2)} 
(z_1,z'_\cap,z'_2),
\end{align}
\end{claim}

\begin{cpf}
Using the definition of $\lambda$, we rewrite \eqref{eq to verify} in the form
\begin{align*}
(\tilde z_1,\tilde z_\cap,\tilde z_2)
& = 
\sum_{\substack{(z_1 ,z_\cap) \in \S_{G_1} \\ (z'_\cap, z'_2) \in \S_{G_2} \\ z_{p_i} = z'_{p_i} \forall i \in [k]}} 
\lambda_{(z_1,z'_\cap,z'_2)} 
(z_1,z'_\cap,z'_2) \\
& =
\sum_{i \in [k+1]} \ 
\sum_{\substack{(z_1 ,z_\cap) \in \S_{G_1} \\ \min\{j \in [k+1] : z_{p_j} = 0\} = i}} \
\sum_{\substack{(z'_\cap, z'_2) \in \S_{G_2} \\ \min\{j \in [k+1] : z'_{p_j} = 0\} = i}} \
\lambda_{(z_1,z'_\cap,z'_2)} (z_1,z'_\cap,z'_2) \\
& = 
\sum_{i \in [k+1]} \ 
\sum_{\substack{(z_1 ,z_\cap) \in \S_{G_1} \\ \min\{j \in [k+1] : z_{p_j} = 0\} = i}} \
\sum_{\substack{(z'_\cap, z'_2) \in \S_{G_2} \\ \min\{j \in [k+1] : z'_{p_j} = 0\} = i}} \
\frac{\mu_{(z_1,z_\cap)} \nu_{(z'_\cap,z'_2)}}{m(i)} (z_1,z'_\cap,z'_2).
\end{align*}
We now verify the obtained inequality, first for components $\tilde z_1$, and then for components $\tilde z_\cap,\tilde z_2$.
We start with components $\tilde z_1$.
Using \cref{claim sum mu nu}, we obtain
\begin{align*}
\tilde z_1 
& = 
\sum_{i \in [k+1]} \ 
\sum_{\substack{(z_1 ,z_\cap) \in \S_{G_1} \\ \min\{j \in [k+1] : z_{p_j} = 0\} = i}} \
\sum_{\substack{(z'_\cap, z'_2) \in \S_{G_2} \\ \min\{j \in [k+1] : z'_{p_j} = 0\} = i}} \
\frac{\mu_{(z_1,z_\cap)} \nu_{(z'_\cap,z'_2)}}{m(i)} z_1\\
& = 
\sum_{i \in [k+1]} 
\frac{1}{m(i)}
\pare{
\sum_{\substack{(z_1 ,z_\cap) \in \S_{G_1} \\ \min\{j \in [k+1] : z_{p_j} = 0\} = i}} 
\mu_{(z_1,z_\cap)}z_1}
\pare{
\sum_{\substack{(z'_\cap, z'_2) \in \S_{G_2} \\ \min\{j \in [k+1] : z'_{p_j} = 0\} = i}} 
\nu_{(z'_\cap,z'_2)}} \\
& = 
\sum_{i \in [k+1]} 
\frac{m(i)}{m(i)} 
\sum_{\substack{(z_1 ,z_\cap) \in \S_{G_1} \\ \min\{j \in [k+1] : z_{p_j} = 0\} = i}} 
\mu_{(z_1,z_\cap)}z_1 \\
& = 
\sum_{(z_1 ,z_\cap) \in \S_{G_1}} 
\mu_{(z_1,z_\cap)} z_1,
\end{align*}
and the resulting equation is implied by \eqref{eq mu conv}.

Next, we consider components $\tilde z_\cap,\tilde z_2$.
Using \cref{claim sum mu nu}, we obtain
\begin{align*}
(\tilde z_\cap,\tilde z_2) 
& = 
\sum_{i \in [k+1]} \ 
\sum_{\substack{(z_1 ,z_\cap) \in \S_{G_1} \\ \min\{j \in [k+1] : z_{p_j} = 0\} = i}} \
\sum_{\substack{(z'_\cap, z'_2) \in \S_{G_2} \\ \min\{j \in [k+1] : z'_{p_j} = 0\} = i}} \
\frac{\mu_{(z_1,z_\cap)} \nu_{(z'_\cap,z'_2)}}{m(i)} (z'_\cap,z'_2)\\
& = 
\sum_{i \in [k+1]} 
\frac{1}{m(i)}
\pare{
\sum_{\substack{(z_1 ,z_\cap) \in \S_{G_1} \\ \min\{j \in [k+1] : z_{p_j} = 0\} = i}} 
\mu_{(z_1,z_\cap)}}
\pare{
\sum_{\substack{(z'_\cap, z'_2) \in \S_{G_2} \\ \min\{j \in [k+1] : z'_{p_j} = 0\} = i}} 
\nu_{(z'_\cap,z'_2)} (z'_\cap,z'_2)} \\
& = 
\sum_{i \in [k+1]} 
\frac{m(i)}{m(i)} 
\sum_{\substack{(z'_\cap, z'_2) \in \S_{G_2} \\ \min\{j \in [k+1] : z'_{p_j} = 0\} = i}} 
\nu_{(z'_\cap,z'_2)} (z'_\cap,z'_2) \\
& = 
\sum_{(z'_\cap, z'_2) \in \S_{G_2}} 
\nu_{(z'_\cap,z'_2)} (z'_\cap,z'_2),
\end{align*}
and the resulting equation is \eqref{eq nu conv}.
\end{cpf}
% The theorem then follows from \cref{claim combination,claim sum one,claim conv comb}
\end{prf}

\begin{figure}[h]
\begin{center}
\includegraphics[width=.95\textwidth]{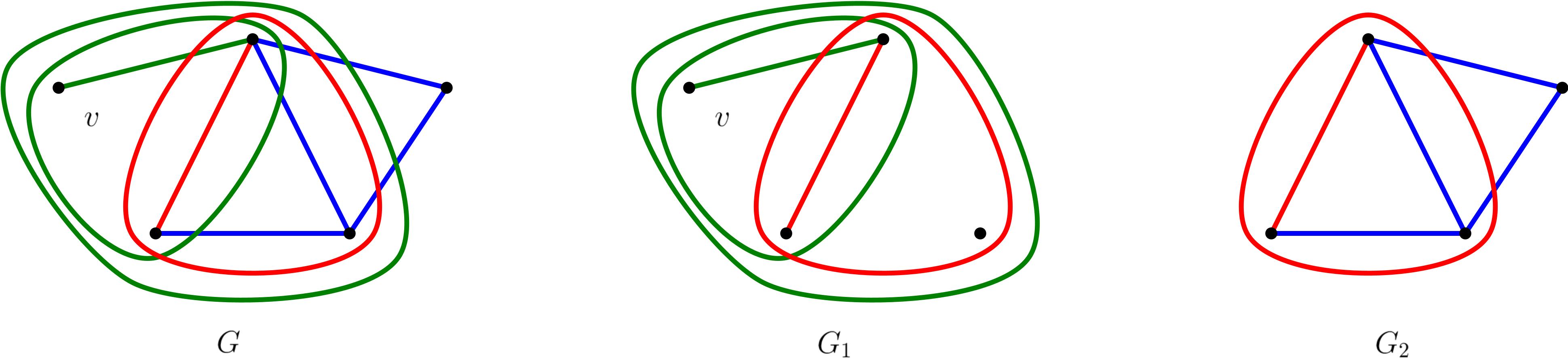}
\end{center}
\caption{An instance of hypergraphs $G, G_1, G_2$ that satisfy the assumptions of \cref{th decomp}: the multilinear set $\S_G$ is decomposable into sets $\S_{G_1}$ and $\S_{G_2}$.}
\label{fig decomp}
\end{figure}

An example of hypergraphs $G, G_1, G_2$ satisfying the assumptions of \cref{th decomp} is given in \cref{fig decomp}.
Theorem~\ref{th decomp} provides a decomposition scheme for multilinear sets whose hypergraphs contain nest points. Moreover, 
the hypergraph $G_1$ defined in the statement of the theorem has a very special structure. In the next section, we characterize the multilinear polytope of $G_1$ using a direct approach. This result together with the decomposition result of Theorem~\ref{th decomp} enables us to obtain a polynomial-size extended formulation for the multilinear polytope of $\beta$-acyclic hypergraphs.

%\newpage

\section{The multilinear polytope of pointed hypergraphs}
\label{sec: pointed}

In this section, we characterize the multilinear polytope for a special type of hypergraphs that serve as the building block for our proposed extended formulations. 
We call  a hypergraph $G=(V,E)$~\emph{pointed} at $v$, if 
\begin{itemize}
\item the edges in $E$ containing $v$ are $e_1 \subset e_2 \subset \dots \subset e_k$, 
\item $V = e_k$,  
\item $E = E_v \cup P_v$, where $E_v := \{e_1, \dots, e_k\}$ and $P_v := \{e_i \setminus \{v\}: i \in [k], \; |e_i\setminus \{v\}| \geq 2\}$. 
\end{itemize}
It then follows that the hypergraph $G_1$ defined in the statement of Theorem~\ref{th decomp} is pointed at $v$. The next theorem provides an explicit description for the multilinear polytope of pointed hypergraphs.

\begin{theorem}
\label{th pointed MPG}
Consider a hypergraph $G = (V, E)$ pointed at $v$.
%Consider a hypergraph $G = (V, E)$ with a nest point $v$ and no isolated nodes. 
%Suppose that $E = E_v \cup P_v$, where $E_v := \{e_i \in E: v \in e_i\}$ and $P_v %:= \{e_i \setminus \{v\}: e_i \in E_v,\; |e_i| \geq 3\}$. 
%Let $k:=|E_v|$ and without loss of generality, let $e_1 \subset e_2 \subset \cdots \subset e_k$. 
Define $p_i = e_i \setminus \{v\}$ for all $i \in [k]$. 
Then $\MP_G$ is defined by the following inequalities: 
\begin{align}
\begin{split}
\label{cvxhull}
    z_u \leq 1 & \qquad \forall u \in V \\
    z_{e_k} \geq 0 \\
    z_{e_k} \leq z_{p_k}  \\
    z_{e_{i+1}} \leq z_{e_i} & \qquad \forall i \in [k-1] \\
    -z_{p_i}+z_{p_{i+1}}+z_{e_i}-z_{e_{i+1}} \leq 0 & \qquad \forall i \in [k-1]\\
    z_{p_{i+1}} \leq z_u & \qquad \forall u \in e_{i+1} \setminus e_{i}, \; \forall i \in [k-1]\\
    \sum_{u \in e_{i+1} \setminus e_{i}}{z_u} + z_{p_{i}}-z_{p_{i+1}} \leq |e_{i+1} \setminus e_{i}| & \qquad \forall i \in [k-1]\\
     z_{e_1} \leq z_{v}\\
    z_{v} + z_{p_1}-z_{e_1} \leq 1\\
     z_{p_1} \leq z_u & \qquad \forall u \in p_1 \\
    \sum_{u \in p_1}{z_u}-z_{p_1} \leq |p_1|-1.
\end{split}
\end{align}
\end{theorem}

\begin{prf}
Denote by $G_0$ (resp. $G_1$) the hypergraph corresponding to the face of $\MP_{G}$ with $z_{v} = 0$ (resp. $z_{v} = 1$). 
We then have: 
$$\MP_G = \conv(\MP_{G_0} \cup \MP_{G_1}).$$ 
It can be checked that both $G_0$ and $G_1$ are $\gamma$-acyclic hypergraphs and hence their multilinear polytopes coincide with their flower relaxations (see Theorem~14 in~\cite{dPKha18SIOPT}).
Denote by $\bar z$ the vector consisting of $z_u$ for all $u \in V \setminus \{v\}$ and $z_e$ for all $e \in P_v$.
It then follows that $\MP_{G_0}$ and $\MP_{G_1}$ are given by:
\begin{align*}
& \MP_{G_0} = \Big\{z \in \R^{|V|+|E|}: z_{v} = 0, \; z_{e_i} =0,  \forall i \in [k], \; \bar z \in \Q \Big\}\\
& \MP_{G_1} = \Big\{z \in \R^{|V|+|E|}: z_{v} = 1, \; z_{e_i} =z_{p_i}, \; \forall i \in [k], \; \bar z \in \Q \Big\},
\end{align*}
where
\begin{align*}
\Q = \Big\{&z \in \R^{|V|+|P_v|-1}: \; z_{u} \leq 1, \; \forall u \in V \setminus \{v\}, \; z_{p_k} \geq 0, \; z_{p_{i+1}} \leq z_{p_i}, \; \forall i \in [k-1], \;  \\
&z_{p_{i+1}} \leq z_u, \; \forall u \in e_{i+1} \setminus e_{i}, \; \forall i \in [k-1], \sum_{u\in e_{i+1} \setminus e_{i}}{z_u}+z_{p_{i}}-z_{p_{i+1}} \leq |e_{i+1} \setminus e_{i}|, \; \forall i \in [k-1]\\
& z_{p_1} \leq z_u, \; \forall u \in p_1, \sum_{u \in p_1}{z_u}-z_{p_1} \leq |p_1|-1 \Big\},
\end{align*}
where the description of $\Q$ follows from Theorem~14 in~\cite{dPKha18SIOPT}.
Using Balas' formulation for the union of polytopes~\cite{Bal98}, it follows that the polytope $\MP_{G}$ is the projection onto the space of the $z$ variables of the polyhedron defined by the following system \eqref{disjc 1}--\eqref{disjc 3}:
\begin{align}
\begin{split}
\label{disjc 1}
\lambda_0 + \lambda_1 = 1, \; \lambda_0 \geq 0, \; \lambda_1 \geq 0\\
z_u = z^0_u + z^1_u & \qquad \forall u \in V \\
z_{p_i} = z^0_{p_i} + z^1_{p_i} & \qquad \forall i \in [k]\\
z_{e_i} = z^0_{e_i} + z^1_{e_i} & \qquad \forall i \in [k]\\
\end{split}
\end{align}
\begin{align}
\begin{split}
\label{disjc 2}
z^0_{v} = 0\\
z^0_{e_i} =0 & \qquad \forall i \in [k]\\
z^0_{u} \leq \lambda_0 & \qquad \forall u \in V \setminus \{v\} \\
z^0_{p_k} \geq 0\\
z^0_{p_{i+1}} \leq z^0_{p_i} & \qquad \forall i \in [k-1]\\
z^0_{p_{i+1}} \leq z^0_u & \qquad \forall u \in e_{i+1} \setminus e_i, \; \forall i \in [k-1]\\
\sum_{u\in e_{i+1} \setminus e_i}{z^0_u}+z^0_{p_{i}}-z^0_{p_{i+1}} \leq |e_{i+1} \setminus e_i| \lambda_0 & \qquad \forall i \in [k-1]\\
z^0_{p_1} \leq z^0_u & \qquad \forall u \in p_1 \\
\sum_{u \in p_1}{z^0_u}-z^0_{p_1} \leq (|p_1|-1) \lambda_0\\
\end{split}
\end{align}
\begin{align}
\begin{split}
\label{disjc 3}
z^1_{v} = \lambda_1 \\
z^1_{e_i} =z^1_{p_i} & \qquad \forall i \in [k]\\
z^1_{u} \leq \lambda_1 & \qquad \forall u \in V \setminus \{v\} \\
z^1_{p_k} \geq 0 \\
z^1_{p_{i+1}} \leq z^1_{p_{i}} & \qquad \forall i \in [k-1]\\
z^1_{p_{i+1}} \leq z^1_u & \qquad \forall u \in e_{i+1} \setminus e_{i}, \; \forall i \in [k-1]\\
\sum_{u\in e_{i+1} \setminus e_i}{z^1_u}+z^1_{p_i}-z^1_{p_{i+1}} \leq |e_{i+1} \setminus e_i| \lambda_1 & \qquad \forall i \in [k-1]\\
z^1_{p_1} \leq z^1_u & \qquad \forall u \in p_1 \\
\sum_{u \in p_1}{z^1_u}-z^1_{p_1} \leq (|p_1|-1) \lambda_1.
\end{split}
\end{align}
In the remainder of this proof, we project out $z^0, z^1, \lambda_0, \lambda_1$ from system~\eqref{disjc 1}--\eqref{disjc 3} and obtain the description of $\MP_{G}$ in the original space. Using $z_{v} = z^0_{v} + z^1_{v}$, $z^0_{v} = 0$, and $z^1_{v} = \lambda_1$, we deduce that $\lambda_0 = 1 -z_{v}$ and $\lambda_1 = z_{v}$. Moreover, from $z_{e_i} = z^0_{e_i} + z^1_{e_i}$, $z^0_{e_i} =0$, and $z^1_{e_i} =z^1_{p_i}$ for all $i \in [k]$, it follows that $z^1_{e_i} =z^1_{p_i}=z_{e_i}$ for all $i \in [k]$, which together with $z_{p_i} = z^0_{p_i} + z^1_{p_i}$ implies that
$ z^0_{p_i} =  z_{p_i}- z_{e_i}$ for all $i \in [k]$. Finally using
$ z_u = z^0_u + z^1_u$ to project out $z^0_u$, for all $u \in V \setminus \{v\}$, the projection of system~\eqref{disjc 1}--\eqref{disjc 3} onto the space $z, z^1_u, u \in V \setminus \{v\}$ is given by:
\begin{align}
0 \leq z_{v} \leq 1 \label{n1}\\
z_{e_k} \geq 0, \; z_{e_k} \leq z_{p_k} \label{n2}\\
z_{e_{i+1}} \leq z_{e_i} & \qquad \forall i \in [k-1] \label{n3}\\
-z_{p_i} +z_{p_{i+1}}+z_{e_i}-z_{e_{i+1}} \leq 0 & \qquad \forall i \in [k-1] \label{n4}
\end{align}
and by 
\begin{align}
\label{remaining}
\begin{split}
z_{u} - z^1_{u}\leq 1-z_{v} & \quad  \forall u \in V \setminus \{v\} \\
z_{p_{i+1}}-z_{e_{i+1}} \leq z_u-z^1_u & \quad   \forall u \in e_{i+1} \setminus e_i, \; \forall i \in [k-1]\\
\sum_{u\in e_{i+1} \setminus e_{i}}{(z_u-z^1_u)}+z_{p_{i}}-z_{e_{i}}-z_{p_{i+1}}+z_{e_{i+1}} \leq |e_{i+1} \setminus e_{i}| (1-z_{v}) & \quad   \forall i \in [k-1]\\
z_{p_1}-z_{e_1} \leq z_u-z^1_u & \quad   \forall u \in p_1 \\
\sum_{u \in p_1}{(z_u-z^1_u)}-z_{p_1}+z_{e_1} \leq (|p_1|-1) (1-z_{v})\\
z^1_{u} \leq z_{v} & \quad   \forall u \in V \setminus \{v\} \\
z_{e_{i+1}} \leq z^1_u & \quad   \forall u \in e_{i+1} \setminus e_{i}, \; \forall i \in [k-1]\\
\sum_{u\in e_{i+1} \setminus e_{i}}{z^1_u}+z_{e_{i}}-z_{e_{i+1}} \leq |e_{i+1} \setminus e_i| z_{v} & \quad   \forall i \in [k-1]\\
z_{e_1} \leq z^1_u & \quad   \forall u \in p_1 \\
\sum_{u \in p_1}{z^1_u}-z_{e_1} \leq (|p_1|-1) z_{v}.
\end{split}
\end{align}
First consider inequalities~\eqref{n1}--\eqref{n4}; besides the redundant inequality $z_{v} \geq 0$, all remaining inequalities are present in system~\eqref{cvxhull}. Hence, to complete the proof, it suffices to project out $z^1_u$, $u \in V \setminus \{v\}$ from system~\eqref{remaining}. 

We start by projecting out variables $z^1_u$, $u \in p_1$ from system~\eqref{remaining}.

\begin{claim}\label{cl1}
Consider all inequalities of system~\eqref{remaining} containing variables $z^1_u$, $u \in p_1$:
\begin{align}
z^1_{u} \leq z_{v} & \qquad \forall u \in p_1 \label{pk1}\\
z_{p_1}-z_{e_1} \leq z_u-z^1_u & \qquad \forall u \in p_1 \label{pk2}\\
\sum_{u \in p_1}{z^1_u}-z_{e_1} \leq (|p_1|-1) z_{v} \label{pk3}\\
z_{u} - z^1_{u}\leq 1-z_{v} & \qquad \forall u \in p_1 \label{pk4}\\
z_{e_1} \leq z^1_u & \qquad \forall u \in p_1 \label{pk5}\\
\sum_{u \in p_1}{(z_u-z^1_u)}-z_{p_1}+z_{e_1} \leq (|p_1|-1) (1-z_{v}). \label{pk6}
\end{align}
Then by projecting out $z^1_u$, $u \in p_1$ from the above system, we obtain 
\begin{align}
    z_u \leq 1 & \qquad \forall u \in p_1\label{p1}\\
    z_{e_1} \leq z_{v} \label{p2}\\
    z_{v} + z_{p_1}-z_{e_1} \leq 1 \label{p3}\\
    z_{p_1} \leq z_u & \qquad \forall u \in p_1 \label{p4}\\
    \sum_{u \in p_1}{z_u}-z_{p_1} \leq |p_1|-1. \label{p5}
\end{align}
%which coincide with inequalities of system~\eqref{cvxhull} containing variables $z_u$, $u \in p_k$, $z_{p_k}$, and $z_{e_k}$.
\end{claim}
\begin{cpf}
 Using Fourier–Motzkin elimination, we project out  $z^1_u$, $u \in p_1$ from inequalities~\eqref{pk1}--\eqref{pk6}. To this end we first consider the following simple cases:
 \begin{enumerate}
     \item Projecting out $z^1_u$, $u \in p_1$ from inequalities~\eqref{pk1} and~\eqref{pk4}, we obtain inequalities~\eqref{p1}.
      \item Projecting out $z^1_u$, $u \in p_1$ from inequalities~\eqref{pk1} and~\eqref{pk5}, we obtain inequality~\eqref{p2}.
      \item Projecting out $z^1_u$, $u \in p_1$ from inequalities~\eqref{pk2} and~\eqref{pk4}, we obtain inequality~\eqref{p3}.
      \item Projecting out $z^1_u$, $u \in p_1$ from inequalities~\eqref{pk2} and~\eqref{pk5}, we obtain inequality~\eqref{p4}. 
\end{enumerate}
To complete the proof, it suffices to project out $z^1_u$, $u \in p_1$ from inequalities~\eqref{pk3} and~\eqref{pk6}. Consider a variable $z^1_{\bar u}$ for some $\bar u \in p_1$. Projecting out this variable from inequalities~\eqref{pk3} and~\eqref{pk6}, we obtain inequality~\eqref{p5}. Hence, to project out $z^1_u$, $u \in p_1$ from inequality~\eqref{pk3} (resp. inequality~\eqref{pk6}), it suffices to consider inequalities~\eqref{pk4} and~\eqref{pk5} (resp. inequalities~\eqref{pk1} and~\eqref{pk2}).  

First consider inequality~\eqref{pk3}; let $p_1 = p^1_1 \cup p^2_1$
such that $p^1_1 \cap p^2_1 = \emptyset$. Projecting out $z^1_u$, $u \in p^1_1$ from inequalities~\eqref{pk3} and~\eqref{pk4},
and projecting out $z^1_u$, $u \in p^2_1$ from inequalities~\eqref{pk3} and~\eqref{pk5}, we obtain:
\begin{equation}\label{red1}
\sum_{u \in p^1_1}{z_u}+(|p^2_1|-1)z_{e_1} \leq (|p^2_1|-1)z_{v}+|p^1_1|. 
\end{equation}
First let $|p^2_1| = 0$; in this case inequality~\eqref{red1} simplifies to
$$
\sum_{u \in e_1}{z_u} - z_{e_1} \leq |e_1|-1,
$$
which is a redundant inequality as it is implied by inequalities~\eqref{p3} and~\eqref{p5}. Now let $|p^2_1| \geq 1$. In this case, inequality~\eqref{red1} is implied by inequalities~\eqref{p1}
and~\eqref{p2} and hence is redundant.

Finally, consider inequality~\eqref{pk6}; let $p_1 = p^1_1 \cup p^2_1$
such that $p^1_1 \cap p^2_1 = \emptyset$. Projecting out $z^1_u$, $u \in p^1_1$ from inequalities~\eqref{pk1} and~\eqref{pk6},
and projecting out $z^1_u$, $u \in p^2_1$ from inequalities~\eqref{pk2} and~\eqref{pk6}, we obtain:
\begin{equation}\label{red2}
\sum_{u \in p^1_1}{z_u}+(|p^2_1|-1) (z_{v} + z_{p_1}-z_{e_1}) \leq (|p_1|-1). 
\end{equation}
First let $|p^2_1| = 0$; in this case inequality~\eqref{red2} simplifies to
$$
\sum_{u \in p_1}{z_u} -z_{v}- z_{p_1} + z_{e_1} \leq |p_1|-1,
$$
which is a redundant inequality as it is implied by inequalities~\eqref{p2} and~\eqref{p5}. Now let $|p^2_1| \geq 1$. In this case, inequality~\eqref{red2} is redundant as it is implied by inequalities~\eqref{p1} and~\eqref{p3}, and this completes the proof.
\end{cpf}

Next, we project out variables $z^1_u$, $u \in e_{i+1} \setminus e_{i}$, $i \in [k-1]$ from system~\eqref{remaining}.

\begin{claim}
\label{cl2}
Let $i \in [k-1]$ and consider all inequalities of system~\eqref{remaining} containing variables $z^1_u$, $u \in e_{i+1} \setminus e_{i}$:
\begin{align}
z^1_{u} \leq z_{v} & \qquad  \forall u \in e_{i+1} \setminus e_{i} \label{pi1}\\
z_{p_{i+1}}-z_{e_{i+1}} \leq z_u-z^1_u & \qquad  \forall u \in e_{i+1} \setminus e_{i} \label{pi2}\\
\sum_{u \in e_{i+1} \setminus e_{i}}{z^1_u}+ z_{e_{i}}-z_{e_{i+1}} \leq |e_{i+1} \setminus e_{i}| z_{v} \label{pi3}\\
z_{u} - z^1_{u}\leq 1-z_{v} & \qquad  \forall u \in e_{i+1} \setminus e_{i} \label{pi4}\\
z_{e_{i+1}} \leq z^1_u & \qquad  \forall u \in e_{i+1} \setminus e_{i} \label{pi5}\\
\sum_{u \in e_{i+1} \setminus e_{i}}{(z_u-z^1_u)}+z_{p_{i}}-z_{e_{i}}-z_{p_{i+1}}+z_{e_{i+1}} \leq (|e_{i+1} \setminus e_{i}|) (1-z_{v}). \label{pi6}
\end{align}
Then by projecting out $z^1_u$, $u \in e_{i+1} \setminus e_i$ from the above system, we obtain a system of inequalities that is implied by  inequalities~\eqref{n3}-\eqref{n4}, inequalities~\eqref{p1}-\eqref{p5} and the following inequalities:
\begin{align}
    z_u \leq 1 & \qquad  \forall u \in e_{i+1} \setminus e_{i}\label{q1}\\
    z_{p_{i+1}} \leq z_{u} & \qquad  \forall u \in e_{i+1} \setminus e_{i}\label{q2}\\
    \sum_{u \in e_{i+1} \setminus e_{i}}{z_u} + z_{p_{i}}-z_{p_{i+1}} \leq |e_{i+1} \setminus e_{i}|. \label{q3}
\end{align}
\end{claim}

\begin{cpf}
Using Fourier–Motzkin elimination, we project out  $z^1_u$, $u \in e_{i+1} \setminus e_{i}$ from inequalities~\eqref{pi1}--\eqref{pi6}.
To this end, it suffices to consider the following cases:
\begin{enumerate}
    \item Projecting out $z^1_u$, $u \in e_{i+1} \setminus e_{i}$ from inequalities~\eqref{pi1} and~\eqref{pi4}, we obtain inequalities~\eqref{q1}.
    \item Projecting out $z^1_u$, $u \in e_{i+1} \setminus e_{i}$ from inequalities~\eqref{pi1} and~\eqref{pi5}, we obtain $z_{e_{i+1}} \leq z_{v}$, which is a redundant inequality as it is implied by
    inequalities~\eqref{n3} and~\eqref{p2}. 
%    $z_{e_i} \leq z_{e_{i+1}}$ and $z_{e_k} \leq z_{v}$.
    \item Projecting out $z^1_u$, $u \in e_{i+1} \setminus e_{i}$ from inequalities~\eqref{pi2} and~\eqref{pi4}, we obtain 
    $z_{v} + z_{p_{i+1}}-z_{e_{i+1}} \leq 1$, whose redundancy follows from inequalities~\eqref{n4} and~\eqref{p3}.
    \item Projecting out $z^1_u$, $u \in e_{i+1} \setminus e_{i}$ from inequalities~\eqref{pi2} and~\eqref{pi5}, we obtain inequalities~\eqref{q2}.
    \item Projecting out $z^1_u$, $u \in e_{i+1} \setminus e_{i}$ from inequalities~\eqref{pi3} and~\eqref{pi6}, we obtain inequalities~\eqref{q3}.
\end{enumerate}
It remains to project out $z^1_u$, $u \in e_{i+1}\setminus e_i$ from inequality~\eqref{pi3} (resp.~\eqref{pi6}) together with inequalities~\eqref{pi4} and~\eqref{pi5} (resp. inequalities~\eqref{pi1} and~\eqref{pi2}).
It can be checked that the resulting inequalities are redundant. We do not include the proof here as it follows from a similar line of arguments to those establishing redundancy of inequalities~\eqref{red1} and~\eqref{red2} in the proof of Claim~\ref{cl1}.
\end{cpf}

Therefore, by inequalities~\eqref{n1}--\eqref{n4}, Claim~\ref{cl1},
and Claim~\ref{cl2}, we conclude that $\MP_{G}$ is defined by system~\eqref{cvxhull}.
\end{prf}

It is interesting to note that, in spite of its simple structure, the constraint matrix of the multilinear polytope of a pointed hypergraph is not totally unimodular. The following example demonstrate this fact. For notational simplicity, in all examples, given a node $v_i$, we write $z_i$ instead of $z_{v_i}$. Similarly, given an edge $\{v_i, v_j, v_k\}$, we write $z_{ijk}$ instead of $z_{\{v_i,v_j, v_k\}}$.

\begin{example}
    Consider $G=(V, E)$ with $V=\{v_1, v_2, v_3, v_4\}$
    and $E=\{\{v_1, v_2\}, \{v_2, v_3, v_4\}, V\}$. It is simple to check that the $G$ is a pointed hypergraph at $v_1$; now consider the following inequalities all of which are present in the description of $\MP_G$:
    \begin{align*}
             z_{234}    \leq z_{3} \\
             z_{234}    \leq  z_{4}\\
        -z_{2}+z_{12}+z_{234}-z_{1234} \leq 0 \\
         z_{2}+z_{3}+z_{4}   -z_{234}  \leq 2.
    \end{align*}
    It can be checked that all above inequalities are facet-defining for $\MP_G$. Now consider the submatrix of these inequalities corresponding to variables $z_2, z_3, z_4, z_{234}$. It can be checked that the determinant of this submatrix equals -2, implying the constraint matrix of $\MP_G$ is not totally unimodular.
\end{example}

We should also remark that one cannot use the concept of balanced matrices to prove the integrality of system~\eqref{cvxhull} (see Theorem 6.13 in~\cite{Cor01b}). In order to use this result, each inequality $a x \leq b$ defining the system should satisfy $b = 1-n(a)$, where $n(a)$ denotes the number of elements in $a$ equal to $-1$. The inequality $-z_{p_i}+z_{p_{i+1}}+z_{e_i}-z_{e_{i+1}} \leq 0$ does not satisfy this assumption as for this inequality we have $n(a)= 2$ and $b = 0$.

%\newpage

\section{Extended formulations and hypergraphs with nest points}
\label{sec: extended}
In this section, we represent the multilinear polytope of hypergraphs with nest points in terms of multilinear polytopes of simpler hypergraphs. As a result, we obtain a  polynomial-size extended formulation for the multilinear polytope of $\beta$-acyclic hypergraphs. To this end, we first introduce expanded hypergraphs, a class of hypergraphs that determine the extended space to which our proposed extended formulations belong.

\subsection{Expanded hypergraphs}
%In this paper we will often consider a special type of hypergraphs.
We say that a hypergraph $G=(V,E)$ is \emph{expanded w.r.t.~$v_1, \dots, v_s$,} if $v_1, \dots, v_s$ is a nest point sequence of $G$, and for every edge $e \in E$, the set $E$ also contains the sets of cardinality at least two among $e \setminus \{v_1\}$, $e \setminus \{v_1, v_2\}$, $\dots$, $e \setminus \{v_1, \dots, v_s\}$. The following three lemmas establish some basic properties of expanded hypergraphs which we will use for our convex hull characterizations:

\begin{lemma}
\label{lem all expanded}
Let $G$ be a hypergraph expanded w.r.t.~$v_1, \dots, v_s$, for $s \ge 1$.
Then, $G - v_1$ is expanded w.r.t.~$v_2, \dots, v_s$.
\end{lemma}

\begin{prf}
Let $G = (V,E)$ and let $G - v_1 = (V',E')$, where $V' := V \setminus \{ v_1 \}$, and $E' := \{ e \setminus \{v_1\} : e \in E, \ |e \setminus \{v_1\}| \ge 2 \}$.
Since $v_1, \dots, v_s$ is a nest point sequence of $G$, then $v_2, \dots, v_s$ is a nest point sequence of $G - v_1$.
Thus, we only need to show that, for every $e \in E'$, the set $E'$ contains the sets of cardinality at least two among $e \setminus \{v_2\}$, $e \setminus \{v_2, v_3\}$, $\dots$, $e \setminus \{v_2, \dots, v_s\}$.

Let $e \in E'$.
We show that we have $e \in E$.
By definition of $G - v_1$, either $e \in E$, or $e \cup \{v_1\} \in E$.
In the first case we are done;
In the second case, since $G$ is expanded, we have $(e \cup \{v_1\}) \setminus \{v_1\} \in E$, thus $e \in E$, and we are done.
Since $e \in E$ and $G$ is expanded, the set $E$ also contains the sets of cardinality at least two among $e \setminus \{v_1\}$, $e \setminus \{v_1, v_2\}$, $\dots$, $e \setminus \{v_1, \dots, v_s\}$.
Since $v_1 \notin e$, these are the sets of cardinality at least two among $e \setminus \{v_2\}$, $e \setminus \{v_2, v_3\}$, $\dots$, $e \setminus \{v_2, \dots, v_s\}$.
By definition of $G - v_1$, the set $E'$ also contains the sets of cardinality at least two among $e \setminus \{v_2\}$, $e \setminus \{v_2, v_3\}$, $\dots$, $e \setminus \{v_2, \dots, v_s\}$.
Hence, $G - v_1$ is expanded w.r.t.~$v_2, \dots, v_s$.
\end{prf}

Let $G = (V,E)$ be a hypergraph and let $v_1, \dots, v_s$ be a nest point sequence of $G$.
The \emph{expansion of $G$ w.r.t.~$v_1, \dots, v_s$} is the hypergraph $G' = (V,E')$, where $E'$ is obtained from $E$ by adding, for each $e \in E$, the sets of cardinality at least two among $e \setminus \{v_1\}$, $e \setminus \{v_1, v_2\}$, $\dots$, $e \setminus \{v_1, \dots, v_s\}$.

\begin{lemma}
\label{lem expansion}
Let $G$ be a hypergraph, let $v_1, \dots, v_s$ be a nest point sequence of $G$,
and let $G'$ be the expansion of $G$ w.r.t.~$v_1, \dots, v_s$.
Then $G'$ is expanded w.r.t.~$v_1, \dots, v_s$.
\end{lemma}

\begin{prf}
Let $G' = (V, E')$.
Clearly, for every edge $e \in E'$, the set $E'$ also contains the sets of cardinality at least two among $e \setminus \{v_1\}$, $e \setminus \{v_1, v_2\}$, $\dots$, $e \setminus \{v_1, \dots, v_s\}$.
Thus, we only need to show that $v_1, \dots, v_s$ is a nest point sequence of $G'$.
Note that, by construction of $G'$, for every $i=1,\dots,s$, the edges of $G - v_1 - \dots - v_{i-1}$ containing $v_i$ coincide with the edges of $G' - v_1 - \dots - v_{i-1}$ containing $v_i$.
Hence, for $i=1,\dots,s$, the fact that $v_i$ is a nest point of $G - v_1 - \dots - v_{i-1}$ implies that $v_i$ is a nest point of $G' - v_1 - \dots - v_{i-1}$.
Thus, $v_1, \dots, v_s$ is a nest point sequence of $G'$.
\end{prf}

\begin{lemma}
\label{lem structure}
Let $G = (V,E)$ be a hypergraph expanded w.r.t.~$v_1, \dots, v_s$.
Let $e \in E$ such that $e \cap \{v_1, \dots, v_s\} \neq \emptyset$, and let $v_i$ be the first node in the sequence $v_1, \dots, v_s$ contained in $e$.
Then, $e \setminus \{v_i\} \in E$, if $|e \setminus \{v_i\}| \ge 2$.
Furthermore, if there exists at least one edge in $E$ strictly contained in $e$ and containing $v_i$, then there exists only one edge $f \in E$ of maximum cardinality, and $f \setminus \{v_i\} \in E$, if $|f| \ge 3$.
\end{lemma}

\begin{prf}
If $|e \setminus \{v_i\}| \ge 2$, we have $e \setminus \{v_i\} \in E$ since $G$ is expanded.

In the rest of the proof, we assume that the set of edges $\bar E := \{g \in E : g \subset e, \ v_i \in g\}$ is nonempty.
We show that there exists one edge $f \in \bar E$ containing all the edges in $\bar E$.

Since $G$ is expanded w.r.t.~$v_1, \dots, v_s$, node $v_i$ is a nest point of $G - v_1 - \cdots - v_{i-1}$.
Hence, the edges of $G - v_1 - \cdots - v_{i-1}$ containing $v_i$ are totally ordered.
%Note that, the edges of $G - v_1 - \cdots - v_{i-1}$ do not contain $v_1,\dots,v_{i-1}$.
Since $G$ is expanded, the edges of $G - v_1 - \cdots - v_{i-1}$ containing $v_i$ coincide with the edges of $G$ containing $v_i$ and not containing $v_1,\dots,v_{i-1}$.
Hence, these edges are totally ordered.
Assume these are edges $e_1 \subset e_2 \subset \cdots \subset e_k$.
Then, $e = e_j$, for $j \ge 2$, and we set $f := e_{j-1}$.

% \bigskip

% , and for every edge $e \in E$, the set $E$ also contains the sets of cardinality at least two among $e \setminus \{v_1\}$, $e \setminus \{v_1, v_2\}$, $\dots$, $e \setminus \{v_1, \dots, v_s\}$.

% To obtain a contradiction, assume that there are two distinct $g_1,g_2 \in E$, strictly contained in $e$, containing $v$, and such that $g_1 \setminus g_2 \neq \emptyset$ and $g_2 \setminus g_1 \neq \emptyset$.
% Let $u_1 \in g_1 \setminus g_2$ and $u_2 \in g_2 \setminus g_1$. 
% Then, $G$ contains the $\beta$-cycle $v,g_1,u_1,e \setminus \{v\},u_2,g_2,v$.
% This contradicts the fact that $G$ is $\beta$-acyclic.
% This concludes the proof that there exists one edge $f \in \bar E$ containing all the edges in $\bar E$.
% This directly implies that $f$ is the only edge in $\bar E$ of maximum cardinality.

Furthermore, since $f$ is contained in $e$ and it contains $v_i$, $v_i$ is the first node in the sequence $v_1, \dots, v_s$ such that $v_i \in f$, thus $f \setminus \{v_i\} \in E$, if $|f| \ge 3$, since $G$ is expanded.
\end{prf}

\subsection{Convex hull characterizations}

In this section, we study the multilinear polytope of expanded hypergraphs.
First we consider the general case in which the hypergraph $G$ is expanded w.r.t. $v_1, \cdots, v_s$ for some $s \geq 1$ and characterize  $\MP_G$ in terms of multilinear polytopes of simpler hypergraphs. Subsequently, we consider the important special case with $s = |V|$ and present a polynomial-size formulation for $\MP_G$. This in turn enables us to obtain a 
polynomial-size extended formulation for the multilinear polytope of $\beta$-acyclic hyerpgraphs. 
Recall that $G'=(V',E')$ is called a~\emph{partial hypergraph} of $G=(V,E)$, if
$V' \subseteq V$ and $E' \subseteq E$.

\begin{theorem}
\label{th general s nodes}
Let $G = (V,E)$ be a hypergraph expanded w.r.t.~$v_1, \dots, v_s$ for some $s \geq 1$.
For each $i \in [s]$, denote by $\tilde G_i$ the partial hypergraph of $G-v_1-v_2-\cdots-v_{i-1}$ pointed at $v_i$; that is, denoting by $E_{v_i}: =\cup_{j \in [k]} {\{e_j\}}$ the set of edges of $G-v_1-v_2-\cdots-v_{i-1}$ containing $v_i$, and letting $e_1 \subset e_2 \subset \cdots \subset e_k$, we have $\tilde G_i = (\tilde V_i,\tilde E_i)$, where $\tilde V_i = e_k$
and $\tilde E_i = E_{v_i}\cup \{e \setminus \{v_i\}: |e\setminus \{v_i\}| \geq 2, e \in E_{v_i}\}$.
Then, $\MP_{G}$ is given by a description of $\MP_{G - v_1 - \cdots - v_s}$ together with a description of $\MP_{\tilde G_i}$ for all $i \in [s]$, where $\MP_{\tilde G_i}$ is characterized in Theorem~\ref{th pointed MPG}.
\end{theorem}

\begin{prf}
The proof is by induction on the number of nest points  $s$ of $G$. In the base case we have $s= 0$. In this case we do not have any pointed hypergraphs $\tilde G_i$
and we have $G - v_1 - \cdots - v_s = G$, hence the statement trivially holds. 

We now show the inductive step. Node $v_1$ is a nest point of $G$; define $G_1 = \tilde G_1$, \ie the partial hypergraph of $G$ pointed at $v_1$, and $G_2 = G - v_1$. Then by Theorem~\ref{th decomp}, the set $\S_G$ is decomposable into $\S_{G_1}$ and $\S_{G_2}$.
That is, $\MP_G$ is defined by inequalities defining $\MP_{\tilde G_1}$ together with those defining $\MP_{G-v_1}$. Since $\tilde G_1$ is a pointed hypergraph at $v_1$, its multilinear polytope $\MP_{\tilde G_1}$ is given by Theorem~\ref{th pointed MPG}. From \cref{lem all expanded}, it follows that $G - v_1$ is expanded w.r.t.~$v_2, \dots, v_{s}$, and it has one fewer nest point than $G$.
Hence, by the induction hypothesis, the polytope $\MP_{G-v_1}$ is given by 
a description of $\MP_{G-v_1- v_2 - \cdots - v_{s}}$ together with
a description of $\MP_{\tilde G_i}$ for all $i \in \{2,\cdots,s\}$, and this completes the proof.
\end{prf}

\begin{theorem}
\label{th2}
Let $G = (V,E)$ be a $\beta$-acyclic hypergraph expanded w.r.t.~$v_1, \dots, v_n$.
For every $e \in E$, we denote by $v(e)$ the first node in the sequence $v_1, \dots, v_n$ contained in $e$, and we define $p(e):=e \setminus \{v(e)\}$.
Define $M := \{e \in E: \exists g \in E, g \subset e, v(e) \in g\}$.
For every $e \in M$, let $f(e) \subset e$ be the edge of maximum cardinality with  $v(e) \in f(e)$ (unique by \cref{lem structure}), and let $f'(e) := f(e) \setminus \{v(e)\}$. 
Finally, denote by $\bar E$ the set of maximal edges of $G$; \ie $\bar E = \{e \in E : \nexists g \in E, \; g \supset e\}$.
Then, $\MP_{G}$ is defined by the following system of linear inequalities:
\begin{align}
0 \le z_{u} \le 1 & \qquad \forall u \in V \label{r1}\\
z_e \ge 0 & \qquad \forall e \in \bar E \label{r2}\\
z_e - z_{p(e)} \le 0 & \qquad \forall e \in E \label{r3}\\
z_e - z_{f(e)} \le 0 & \qquad \forall e \in M\label{r4}\\
- z_{f'(e)} + z_{p(e)} + z_{f(e)} - z_e \le 0 & \qquad \forall e \in M\label{r5}\\
z_{v(e)} + z_{p(e)} - z_e \le 1 & \qquad \forall e \in E \setminus M \label{r6}\\
z_e - z_{v(e)} \le 0 & \qquad \forall e \in E \setminus M\label{r7}.
\end{align}
\end{theorem}

\begin{prf}
First note that the all variables that appear in system~\eqref{r1}--\eqref{r7} are present in $\S_G$ due to \cref{lem structure}.
The proof is by induction on $|V|$.
In the base case we have $V = \{u\}$ and $E = \emptyset$.
Clearly, $\MP_{G}$ is then given by $0 \le z_u \le 1$.

We now show the inductive step.
Node $v_1$ is a nest point of $G$, and for ease of notation we set $v := v_1$.
Without loss of generality, assume that $v$ is not an isolated node.
Let $e_1 \subset e_2 \subset \cdots \subset e_k$ for some $k \geq 1$, be the edges of $G$ containing $v$.
For each $i \in [k]$, let $p_i := e_i \setminus \{v\}$,
let $E_v := \{e_1,\dots,e_k\}$, and let $P_v := \{p \in \{p_1,\dots,p_k\} : |p| \ge 2\}$.
Define $G_1 := (e_k, E_v \cup P_v)$ and $G_2 := G - v$.
Since $G$ is expanded w.r.t. $v_1, \cdots, v_n$, we have $G = G_1 \cup G_2$.
Then from \cref{th decomp}, it follows that the set $\S_G$ is decomposable into $\S_{G_1}$ and $\S_{G_2}$. That is, $\MP_G$ is defined by inequalities defining $\MP_{G_1}$ together with those defining $\MP_{G_2}$.

By~\cref{th pointed MPG}, the polytope $\MP_{G_1}$ is given by:
\begin{align}
    0 \leq z_u \leq 1 & \qquad \forall u \in e_k \label{h1}\\
    z_{e_k} \geq 0 \label{h2}\\
%    z_{e_k} -z_{p_k} \leq 0 \label{h3}\\
    z_{e_i} -z_{p_i} \leq 0 & \qquad \forall i \in [k] \label{h3}\\
    z_{e_{i+1}} -z_{e_i} \leq 0 & \qquad \forall i \in [k-1] \label{h4}\\
    -z_{p_i}+z_{p_{i+1}}+z_{e_i}-z_{e_{i+1}} \leq 0 & \qquad \forall i \in [k-1] \label{h5}\\
    z_{v} + z_{p_1}-z_{e_1} \leq 1 \label{h6}\\
    z_{e_1} -  z_{v}\leq 0 \label{h7}\\
    z_{p_{i+1}} \leq z_u & \qquad \forall u \in e_{i+1} \setminus e_{i}, \; \forall i \in [k-1] \label{h8}\\
    \sum_{u \in e_{i+1} \setminus e_{i}}{z_u} + z_{p_{i}}-z_{p_{i+1}} \leq |e_{i+1} \setminus e_{i}| & \qquad \forall i \in [k-1]\label{h9}\\
    z_{p_1} \leq z_u & \qquad \forall u \in p_1 \label{h10}\\
    \sum_{u \in p_1}{z_u}-z_{p_1} \leq |p_1|-1. \label{h11}
\end{align}
We should remark that the valid inequalities $z_u \geq 0$ for all $u \in e_k$ and $z_{e_i} \leq z_{p_i}$, for all $i \in [k-1]$ are not present in the description of $\MP_{G_1}$
as given by~\cref{th pointed MPG}, and hence are redundant. However, we includ them in the above system as they simplify the proof.  
From \cref{lem all expanded}, it follows that $G - v_1$ is a $\beta$-acyclic hypergraph expanded w.r.t.~$v_2, \dots, v_{n}$, and it has one fewer node than $G$.
Hence, by the induction hypothesis, the polytope $\MP_{G_2}$ is given by
\begin{align}
0 \le z_{u} \le 1 \qquad &\forall u \in V \setminus \{v\}\label{g1}\\
z_e \ge 0 \qquad &\forall e \in \bar E \setminus \{e_k\}\label{g2}\\
z_e - z_{p(e)} \le 0 \qquad &\forall e \in E \setminus E_v\label{g3}\\
z_e - z_{f(e)} \le 0 \qquad &\forall e \in M \setminus \cup_{i=2}^k {\{e_i\}}\label{g4}\\
- z_{f'(e)} + z_{p(e)} + z_{f(e)} - z_e \le 0 \qquad &\forall e \in M \setminus \cup_{i=2}^k {\{e_i\}}\label{g5}\\
z_{v(e)} + z_{p(e)} - z_e \le 1 \qquad &\forall e \in E \setminus (M \cup \{e_1\}) \label{g6}\\
z_e - z_{v(e)} \le 0  \qquad &\forall e \in E \setminus (M \cup \{e_1\}). \label{g7}
\end{align}

Hence, to complete the proof it suffices to show that combining inequalities~\eqref{h1}--\eqref{h11} and inequalities~\eqref{g1}--\eqref{g7}, we obtain inequalities~\eqref{r1}--\eqref{r7}. We start by making the following observations:
\begin{enumerate}
    \item Inequalities~\eqref{h1} and~\eqref{g1} are equivalent to inequalities~\eqref{r1}.
    \item Inequalities~\eqref{h2} and~\eqref{g2} are equivalent to inequalities~\eqref{r2}.
    \item Inequalities~\eqref{h3} and~\eqref{g3} are equivalent to inequalities~\eqref{r3}.
    \item Noting that $f(e_{i+1}) = e_i$, for all $i \in [k-1]$, it follows that
inequalities~\eqref{h4} and~\eqref{g4} are equivalent to inequalities~\eqref{r4}.
\item Noting that $f'(e_{i+1}) = p_i$, for all $i \in [k-1]$, it follows that
inequalities~\eqref{h5} and~\eqref{g5} are equivalent to inequalities~\eqref{r5}.
\item Inequalities~\eqref{h6} and~\eqref{g6} are equivalent to inequalities~\eqref{r6}.
\item inequalities~\eqref{h7} and~\eqref{g7} are equivalent to inequalities~\eqref{r7}.
\end{enumerate}

Hence it remains to show that inequalities~\eqref{h8}--\eqref{h11} are implied by inequalities~\eqref{r1}--\eqref{r7}. First, let us consider inequalities~\eqref{h8};
namely,
\begin{equation}\label{redd1}
z_{p_{i+1}} \leq z_u \qquad u \in e_{i+1} \setminus e_{i}, \; i \in [k-1].
\end{equation}
Fix $i \in [k-1]$ and fix $u \in p_{i+1}$; note that since $G$ is expanded w.r.t. $v$, we have $p_{i+1} \in E$. For notational simplicity let $q_0 := p_{i+1}$.
Consider the two elements $p(q_0), v(q_0)$, if $q_0 \in E \setminus M$ (resp. $p(q_0), f(q_0)$, if $q_0 \in  M$), as defined in system~\eqref{g1}-~\eqref{g7}. Let $q_1$ be one of these two elements such that $u \in q_1$. Clearly the inequality $z_{q_0} \leq z_{q_1}$
is present in system~\eqref{g1}-~\eqref{g7}. If $q_1 = u$, then we are done. Otherwise, as before we consider 
$p(q_1), v(q_1)$, if $q_1 \in E \setminus M$ (resp. resp. $p(q_1), f(q_1)$, if $q_1 \in M$), and let $q_2$ be one of these two elements such that $u \in q_2$. Again the inequality $z_{q_1} \leq z_{q_2}$ is present in system~\eqref{g1}-~\eqref{g7}. We continue this recursion until $q_k = u$. It then follows that the collection of inequalities $z_{p_{i+1}}=z_{q_0} \leq z_{q_1}$, $z_{q_1} \leq z_{q_2}, \cdots, z_{q_{k-1}} \leq z_{q_k}= z_u$ imply
inequality~\eqref{redd1}.
%
%Consider the following inequalities of system~\eqref{r1}--\eqref{r7}:
%\begin{align*}
%& z_e \leq z_{p(e)} \qquad  \forall e \in E: e \subseteq p_{i+1}\\
%& z_e \leq z_{v(e)} \qquad  \forall e \in M: e \subseteq p_{i+1}\\
%& z_e \leq z_{f(e)} \qquad  \forall e \in E \setminus M: e \subseteq p_{i+1}.
%\end{align*}
%It then follows that the above inequalities imply $z_{p_{i+1}} \leq z_u$ for all $u \in p_{i+1}$, which in turn imply inequalities~\eqref{redd1} are redundant.
The redundancy of inequalities~\eqref{h10} follows from a similar line of arguments.

Next, let us consider inequalities~\eqref{h9}, which can be written as:
\begin{equation}\label{redd2}
\sum_{u \in p_{i+1} \setminus p_{i}}{z_u} + z_{p_{i}}-z_{p_{i+1}} \leq |p_{i+1} \setminus p_{i}| \qquad  i \in [k-1],
\end{equation}
where as before for any $i \in [k-1]$, we assume $|p_{i+1}| \geq 2$ and hence by construction
$p_{i+1} \in E$. In the following we show that these inequalities are implied by system~\eqref{r1}--\eqref{r7}. To this end, we prove a more general statement, \ie for any $q_1 \subset q_2$ with $q_1 \in V \cup E$ and $q_2 \in E$, we show that the inequality
\begin{equation}\label{rednew}
\sum_{u \in q_2 \setminus q_1}{z_u}+z_{q_1}-z_{q_2} \leq |q_2 \setminus q_1|,
\end{equation}
is implied by the following inequalities of system~\eqref{r1}--\eqref{r7}:
\begin{align}
-  z_{f'(e)} + z_{p(e)} + z_{f(e)} - z_e \le 0 & \qquad \forall e \in M: e \subseteq q_2, \; e \not\subseteq q_1\label{tt1}\\
z_{v(e)} + z_{p(e)} - z_e \le 1 & \qquad \forall e \in E \setminus M: e \subseteq q_2, \; e \not\subseteq q_1.\label{tt2}
\end{align}
Then setting $q_1=p_i$ and $q_2 = p_{i+1}$ completes the proof.
The proof is by induction on the number of nodes in $q_2$. In the base case we have $q_2=\{u_1, u_2\}$ for which inequality~\eqref{rednew} simplifies to
$z_{u_1} + z_{u_2}-z_{q_2} \leq 1$, which is present among inequalities~\eqref{tt2} since from $|q_2|=2$, it follows that $q_2 \in E \setminus M$. 

We now proceed with the inductive step. Let $|q_2|=k$ for some $k \geq 3$. Two cases arise:
\begin{itemize}
\item [(i)] $q_1 \subseteq p(q_2)$: in this case the following inequality is present among (if $q_1 \in E\setminus M$) or is implied by (if $q_1 \in M$) inequalities~\eqref{tt1} and~\eqref{tt2}:
\begin{equation}\label{redx1}
z_{v(q_2)}+ z_{p(q_2)}-z_{q_2} \leq 1.
\end{equation}
Since by assumption $q_1 \subseteq p(q_2)$ and $|p(q_2)|=k-1$, by the induction hypothesis, the inequality
\begin{equation}\label{redx2}
\sum_{u \in p(q_2) \setminus q_1}{z_u} + z_{q_1}-z_{p(q_2)} \leq |p(q_2) \setminus q_1|,
\end{equation}
is implied by inequalities
\begin{align*}
-  z_{f'(e)} + z_{p(e)} + z_{f(e)} - z_e \le 0 & \qquad \forall e \in M: e \subseteq p(q_2), \; e \not\subseteq q_1\\
z_{v(e)} + z_{p(e)} - z_e \le 1 & \qquad \forall e \in E \setminus M: e \subseteq p(q_2), \; e \not\subseteq q_1,
\end{align*}
which are in turn present among inequalities~\eqref{tt1} and~\eqref{tt2}, since $q_1 \subseteq p(q_2) \subseteq q_2$. Summing up inequalities~\eqref{redx1} and~\eqref{redx2} we obtain inequality~\eqref{rednew}.
\item [(ii)] $q_1 \not\subseteq p(q_2)$: in this case, we must have $q_2 \in M$. Hence the following inequality is present among inequalities~\eqref{tt1}:
\begin{equation}\label{redx3}
-z_{f'(q_2)} + z_{f(q_2)} + z_{p(q_2)}-z_{q_2} \leq 0.
\end{equation}
By Lemma~\ref{lem structure}, we have $q_1 \subseteq f(p_2)$. Then using a similar line of arguments to those in case~(i) above, we conclude that the following are implied by inequalities~\eqref{tt1} and~\eqref{tt2}:
\begin{align}
\begin{split}\label{redx4}
& \sum_{u \in p(q_2) \setminus f'(q_2)}{z_u}+z_{f'(q_2)}-z_{p(q_2)} \leq |p(q_2) \setminus f'(q_2)|\\
& \sum_{u \in f(q_2) \setminus q_1)}{z_u}+z_{q_1}-z_{f(q_2)} \leq |f(q_2) \setminus q_1|.
\end{split}
\end{align}
Summing up inequalities~\eqref{redx3} and~\eqref{redx4} we obtain inequality~\eqref{rednew} implying it is redundant.
\end{itemize}

%We claim that summing up the above inequalities, we obtain inequalities~\eqref{redd2}, implying their redundancy. To see this,
The redundancy of inequalities~\eqref{h11} then immediately follows by setting $q_1 = v$ for some $v \in p_1$ and $q_2 = p_1$.
\end{prf}

Let $G=(V,E)$ be a $\beta$-acyclic hypergraph expanded with respect to $v_1, \cdots, v_n$, where $n = |V|$.
By \cref{th2}, $\MP_G$ is given by system~\eqref{r1}--\eqref{r7}.
We should remark that, in spite of its simplicity, the constraint matrix of $\MP_{G}$ is \emph{not} totally unimodular. The following example demonstrates this fact.

\begin{example}
Consider the hypergraph $G=(V, E)$ with $V= \{v_1, v_2, v_3, v_4\}$ and 
$$E= \{\{v_1,v_2\}, \{v_2, v_3\}, \{v_3, v_4\}, \{v_1, v_2, v_3\}, V\}.$$ 
It is simple to check that $G$ is $\beta$-acyclic and is expanded with respect to $v_1, v_2, v_3, v_4$.
By Theorem~\ref{th2}, $\MP_G$ contains the following inequalities:
\begin{align*}
z_1+z_2      -z_{12}                          & \leq 1\\
z_2+z_3      -z_{23}                     & \leq 1\\
                z_{123}-z_{12}        & \leq 0\\
-z_3   +z_{34}+z_{123}-z_{1234} & \leq 0.
\end{align*}
It can be checked that all above inequalities are facet-defining.
The determinant of the submatrix corresponding to variables $z_2,z_3,z_{12},z_{123}$ equals $-2$, implying the constraint matrix of $\MP_{G}$ is not totally unimodular.
\end{example}

Moreover, one cannot use the concept of balanced matrices to prove the integrality of system~\eqref{r1}--\eqref{r7} (see Theorem 6.13 in~\cite{Cor01b}). In order to use this result, each inequality $a x \leq b$ defining the system should satisfy $b = 1-n(a)$, where $n(a)$ denotes the number of elements in $a$ equal to $-1$. Clearly, the inequality $-z_{f'(e)}+z_{p(e)}+z_{f(e)}-z_{e} \leq 0$ does not satisfy this assumption.

\medskip
From Theorems~\ref{th general s nodes} and~\ref{th2}, we directly obtain the following results on extended formulations of multilinear polytopes:

\Second*
%\begin{theorem}
%\label{th general sequence}
%Let $G = (V,E)$ be a hypergraph of rank $r$, and let %$v_1,\dots,v_s$ be a nest point sequence of $G$.
%Then an extended formulation of $\MP_G$ is given by a %description of $\MP_{G - v_1 - \cdots - v_s}$, together %with a system of at most $|V|+(4s+2) |E| + 4 (r-2)s^2$ linear inequalities, including at most $(r-2) s$ extended variables.
%The system is characterized in \cref{th general s nodes}.
% Then, we can construct in polynomial time a hypergraph $G^\circ$, obtained from $G$ by adding at most $s|E|$ edges, so that $\MP_{G^\circ}$ is given by a description of $\MP_{G - v_1 - \cdots - v_s}$ and the following system of linear inequalities:
% A description of $\MP_G$ can then be obtained from the description of $\MP_{G^\circ}$ by projecting out the additional edges introduced.
%\end{theorem}

\begin{prf}
Let $G' = (V, E')$ be the expansion of $G$ w.r.t.~$v_1, \dots, v_s$.
From \cref{lem expansion}, it follows that $G'$ is expanded w.r.t.~$v_1, \dots, v_s$.
We then apply \cref{th general s nodes} to $G'$, and observe that $G' - v_1 - \cdots - v_s = G - v_1 - \cdots - v_s$. The total number of inequalities associated with multilinear polytopes of pointed partial hypergraphs of $G'$ at $v_i$, $i \in [s]$
is upper bounded by $|V|+2 |\bar E|+ 4rs$, where $\bar E$ denotes the set of maximal edges of $G'$. To see this, consider system~\eqref{cvxhull} defining the convex hull of a pointed hypergraph. First, note that we have a total number of $|V|$ inequalities of the form $z_u \leq 1$. The total number of nonredundant inequalities of the form  $z_{e_k} \geq 0$ and $z_{e_k} \leq z_{p_k}$ is $2|\bar E|$. Moreover, the total number of inequalities of the form  $z_{e_{i+1}} \leq z_{e_i}$,
$i \in [k-1]$ and $z_{e_1} \leq z_{v}$ is upper bounded by $rs$. Similarly the total number of inequalities of the form $-z_{p_i}+z_{p_{i+1}}+z_{e_i}-z_{e_{i+1}} \leq 0$, $i \in [k-1]$ and $z_{v} + z_{p_1}-z_{e_1} \leq 1$
is upper bounded by $rs$. Also, the total number of inequalities of the form $z_{p_{i+1}} \leq z_u$, $u \in e_{i+1} \setminus e_{i}$, $i \in [k-1]$, $z_{p_1} \leq z_u$, $u \in p_1$ is upper bounded by $rs$. Finally, the total number of inequalities of the form $\sum_{u \in e_{i+1} \setminus e_{i}}{z_u} + z_{p_{i}}-z_{p_{i+1}} \leq |e_{i+1} \setminus e_{i}|$, $i \in [k-1]$ and $\sum_{u \in p_1}{z_u}-z_{p_1} \leq |p_1|-1$ is upper bounded by $rs$. An upper bound on the number of corresponding  linear inequalities can then be obtained using $|\bar E| \leq |E|$.
Finally, notice that in a rank $r$ hyperpgraph each nest point is present in at most $r-1$ edges implying that the number of extended variables does not exceed $(r-2)s$.

\end{prf}

We are now ready to prove the main result of this paper, \cref{th main}, which we recall below.

\First*
% \begin{corollary}\label{cor 2}
% Let $G = (V,E)$ be a $\beta$-acyclic hypergraph of rank $r$.
% Then there exists an extended formulation of $\MP_G$ comprising of at most 
% $(3r-4)|V|+4|E|$ inequalities, with at most $(r-2) |V|$ extended variables.
% The system is explicitly given in~\cref{th2}.
% \end{corollary}

\begin{prf}
Since $G$ is $\beta$-acyclic, by~\cref{th beta iff}, it has a nest point sequence of length $|V|$, say $v_1, \dots, v_n$.
Let $G' = (V, E')$ be the expansion of $G$ w.r.t.~$v_1, \dots, v_n$.
From \cref{lem expansion}, $G'$ is expanded w.r.t.~$v_1, \dots, v_n$.
We then apply \cref{th2} to $G'$. System~\eqref{r1}--\eqref{r7} consists of $2|V|+ 3|E'|+|\bar E|$ inequalities, where $\bar E$ denotes the set of maximal edges of $G'$. The result then follows using the fact that $|E'| \leq (r-2)|V| +|E|$ and $\bar E \leq E$.
\end{prf}

We remark that \cref{int2} allows us to obtain an extended formulation for the multilinear polytope of certain hypergraphs that are not $\beta$-acyclic.
This happens precisely when a description of $\MP_{G - v_1 - \cdots - v_s}$ is available.
The following example demonstrates this fact.

\begin{figure}[h]
\begin{center}
\includegraphics[width=.35\textwidth]{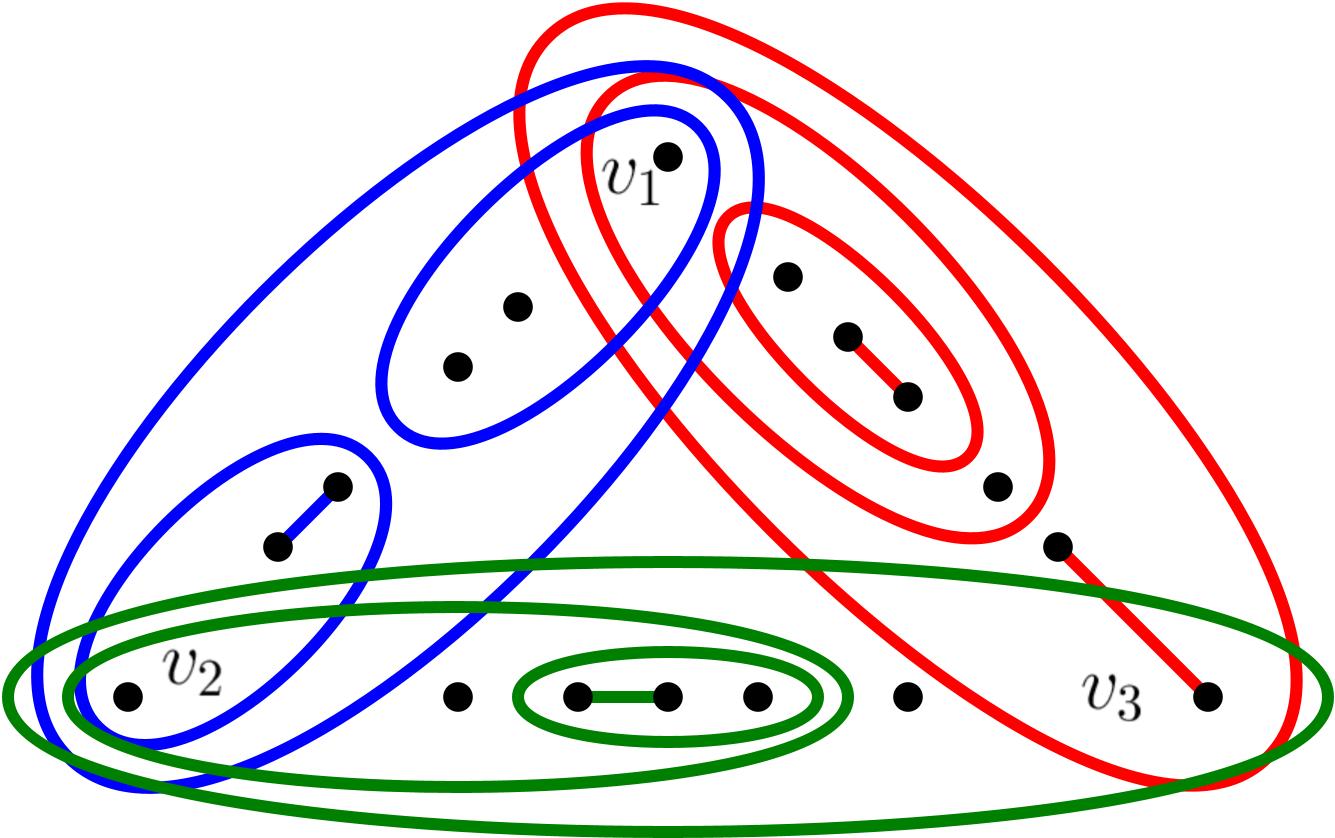}
\end{center}
\caption{Illustration of the hypergraph considered in \cref{ex hypertriangle}.}
\label{fig hypertriangle}
\end{figure}

\begin{example}
\label{ex hypertriangle}
Consider the hypergraph $G$ depicted in \cref{fig hypertriangle}. 
A nest point sequence of $G$ is given by all nodes of $G$ (in any order), except for $v_1,v_2,v_3$.
The hypergraph $G' = (V',E')$ obtained from $G$ by removing all nodes except for $v_1,v_2,v_3$ is a ``triangle'', i.e. $V' = \{v_1,v_2,v_3\}$ and $E' = \{\{v_1, v_2\}, \{v_2, v_3\}, \{v_3, v_1\}\}$.
It is well known that $\MP_{G'}$ is obtained by adding triangle inequalities to the standard linearization  $\MP^{\rm LP}_{G'}$~\cite{Pad89}.
\cref{int2} then gives an extended formulation of $\MP_{G}$.
\end{example}

%\section{The multilinear polytope of $\beta$-acyclic hypergraphs}
%\label{sec projection}

%In this section, we obtain an explicit description for the multilinear polytope of $\beta$-acyclic hypergraphs in the original space. To this end, we first present a natural generalization of running intersection inequalities, a class of valid inequalities for multilinear sets introduced in~\cite{dPKha21MOR}.

\section{The original space}
\label{sec:original}

In Section~\ref{sec: extended}, we presented a polynomial-size extended formulation for the multilinear polytope of $\beta$-acyclic hypergraphs. It is often desirable to obtain an explicit description for the multilinear polytope in the original space. To this end, one can employ Fourier-Motzkin elimination to project out the extended variables from system~\eqref{r1}--\eqref{r7}. 
In~\cite{dPKha18SIOPT}, the authors show that in the original space, the multilinear polytope of a $\gamma$-acyclic hypergraph $G = (V,E)$ contains exponentially many facet-defining inequalities (as a function of $|V|,|E|$), in general.  As $\beta$-acyclicity subsumes $\gamma$-acyclicity, this result implies that the multilinear polytope of a $\beta$-acyclic hypergraph contains
exponentially many facet-defining inequalities, in general. 

%On the other hand, while the multilinear polytope of a $\gamma$-acyclic hypergraph in the original space consists of sparse facet-defining inequalities, as we detail next, the multilinear polytope of a $\beta$-acyclic hypergraph may contain very dense facets in the original space.

From a computational perspective, sparsity is key to the effectiveness of cutting planes in a branch-and-cut framework. Indeed, all existing families of cutting planes for multilinear sets, such as flower inequalities~\cite{dPKha18SIOPT} and running intersection inequalities~\cite{dPKha21MOR} are sparse. Namely, for a rank $r$ hypergraph, flower inequalities contain at most $\frac{r}{2}$ nonzero coefficients, and running intersection inequalities contain at most $2(r-1)$ nonzero coefficients.  When added to the standard linearization, flower inequalities characterize the multilinear polytope of $\gamma$-acyclic hypergraphs~\cite{dPKha18SIOPT}, and running intersection inequalities characterize the multilinear polytope of kite-free $\beta$-acyclic hypergraphs~\cite{dPKha21MOR}.
However, as we detail in the following, the multilinear polytope of a $\beta$-acyclic hypergraph $G=(V,E)$ may contain very dense facets, in general. That is, inequalities containing as many as $\theta(|E|)$ nonzero coefficients. This is significant, as almost for all multilinear sets appearing in nonconvex problems, we have $r \ll |E |$.

\subsection{The multilinear polytope of beta-acyclic hypergraphs with dense facets}

In the following, we present a family of $\beta$-acyclic hypergraphs $G=(V, E)$ whose multilinear polytope contains facet-defining inequalities with $|E|$ non-zero coefficients.

\begin{figure}[h]
\includegraphics[width=.95\textwidth]{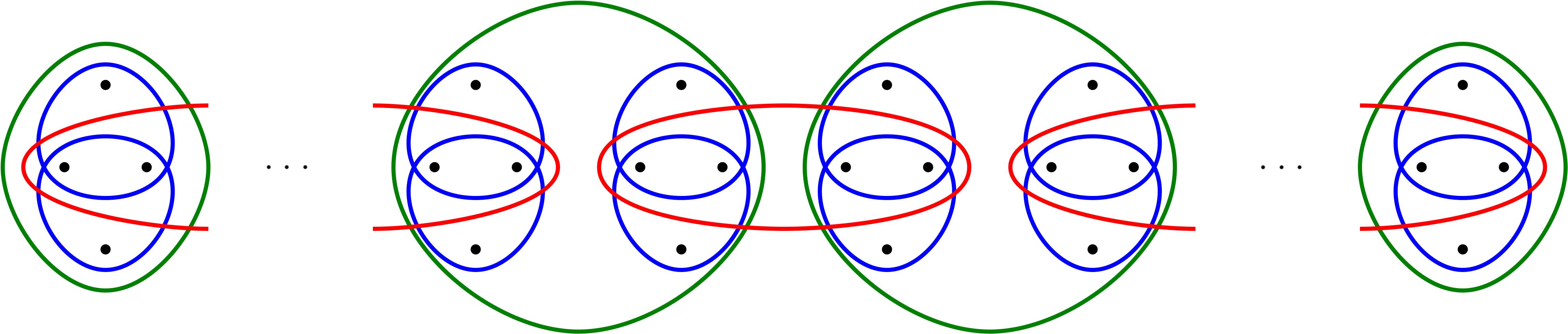}
\caption{Illustration of the family of hypergraphs considered in \cref{prop dense}.}
\label{fig path}
\end{figure}

\begin{proposition}
\label{prop dense}
Let $n \geq 2$ and consider the $\beta$-acyclic hypergraph $G=(V,E)$
with 
$$
V= \bigcup_{i \in [n]}{V^i}, \qquad E= H \cup \bigcup_{i\in [n]}{E^i},
$$
where $V^1 = \{v^1_3, v^1_4, v^1_7, v^1_8\}$, $V^i = \{v^i_1,\cdots, v^i_8\}$ for all $i \in [n-1] \setminus \{1\}$, $V^n = \{v^n_1, v^n_2, v^n_5, v^n_6\}$, 
\begin{align*} 
& H = \Big\{\{v^i_3, v^i_4, v^{i+1}_1, v^{i+1}_2\}, \; i \in [n-1]\Big\}\\
&E^1 = \Big\{\{v^1_3,v^1_4,v^1_7\}, \{v^1_3,v^1_4,v^1_8\}, V^1\Big\} \\
& E^i =\Big\{\{v^i_1,v^i_2,v^i_5\}, \{v^i_1,v^i_2,v^i_6\}, 
\{v^i_3,v^i_4,v^i_7\}, \{v^i_3,v^i_4,v^i_8\}, V^i\Big\} \qquad \forall i \in [n-1] \setminus \{1\}\\
& E^n = \Big\{\{v^n_1,v^n_2,v^n_5\}, \{v^n_1,v^n_2,v^n_6\}, V^n\Big\}.
\end{align*}
See \cref{fig path}. Then the following inequality containing $|E|$ nonzero coefficients defines a facet of $\MP_G$:
\begin{equation}\label{densefacet}
-\sum_{i \in [n]}{z_{V^i}}-\sum_{e \in H}{z_e}+\sum_{i \in [n]}{\sum_{e \in E^i \setminus \{V^i\}} {z_e}} \leq 2 n -3.
\end{equation}
\end{proposition}
\begin{prf}
For notational simplicity in the following, we define $h^i := \{v^i_3, v^i_4, v^{i+1}_1, v^{i+1}_2\}$, for all $i \in [n-1]$, 
$e^i_1 := \{v^i_1,v^i_2,v^i_5\}$, and
$e^i_2 := \{v^i_1,v^i_2,v^i_6\}$ for all $i \in [n] \setminus \{1\}$.
Moreover, we define $e^i_3 := \{v^i_3,v^i_4,v^i_7\}$, and
$e^i_4 := \{v^i_3,v^i_4,v^i_8\}$ for all $i \in [n-1]$.

We start by proving the validity of inequality~\eqref{densefacet} for $\MP_G$. 
First, we construct the hypergraph $G'=(V, E \cup E')$, where 
$$
E':=\Big\{f^i:=\{v^i_1, v^i_2\}, i\in [n] \setminus \{1\}\Big\} \cup \Big\{g^i:=\{v^i_3,v^i_4\}, i \in [n-1]\Big\}.
$$
The following inequalities are all extended running intersection inequalities and hence are valid for $\MP_{G'}$ (see Section~\ref{extendedrunning} for the definition of extended running intersection inequalities):
\begin{align*}
-&z_{g^1}+z_{e^1_3}+z_{e^1_4}-z_{V^1} \leq 0\\
-&z_{f^n}+z_{e^n_1}+z_{e^n_2}-z_{V^n} \leq 0\\
&z_{g^i}+z_{f^{i+1}}-z_{h^i} \leq 1, \quad \forall i \in [n-1]\\
-&z_{f^i}-z_{g^i}+z_{e^i_1}+z_{e^i_2}+z_{e^i_3}+z_{e^i_4}
-z_{V^i} \leq 1, \quad \forall i \in [n-1]\setminus\{1\}.
\end{align*}
Summing up the above inequalities we obtain inequality~\eqref{densefacet} implying its validity for $\MP_{G'}$. Since $\MP_{G'} \subset \MP_{G}$, we conclude that inequality~\eqref{densefacet} is valid for $\MP_G$ as well.

We now show that inequality~\eqref{densefacet} defines a facet of $\MP_G$.
Consider a nontrivial valid inequality
$a z \leq \alpha$ for $\MP_G$ that is
satisfied tightly by any point in $\S_G$ satisfying inequality~\eqref{densefacet} tightly.
In the following, we show that the two inequalities~\eqref{densefacet} and $a z \leq \alpha$ coincide up to a
positive scaling, which by full dimensionality of $\MP_G$ (see Proposition~1 in~\cite{dPKha17MOR}) implies that inequality~\eqref{densefacet} is defines a facet of $\MP_G$.

First consider a point $z^1 \in \S_G$ with $z^1_{v^i_3}=  z^1_{v^i_4} = z^1_{v^i_7} = 1$ for all $i \in [n-1]$, $z^1_{v^i_8}=1$ for all $i \in [n-1]\setminus \{1\}$, and for every other $v \in V$, we have $z^1_v = 0$. It can be checked that inequality~\eqref{densefacet} is satisfied tightly at this point. Now consider a second tight point $z^2 \in \S_G$ whose components are equal to $z^1$ except for one component $z^2_{v^j_1} = 1$ for some $j \in [n] \setminus \{1\}$. Substituting these two tight points in $a z = \alpha$, yields $a_{v^j_1} = 0$. Using a similar line of arguments, we obtain:
\begin{equation}\label{ti1}
a_v = 0 \qquad \forall v \in V.
\end{equation}
Let us again consider the tight point $z^1 \in \S_G$ defined above. Construct another tight point $z^3 \in \S_G$ with $z^3_{v^i_3}=  z^3_{v^i_4} = z^3_{v^i_7} = 1$ for all $i \in [n-1]$, $z^3_{v^i_8}=1$ for all $i \in [n-1] \setminus \{1\}$, $z^3_{v^j_1}=  z^3_{v^j_2} = z^3_{v^j_5} = 1 $, for some $j \in [n]\setminus \{1\}$, and for every other $v \in V$, we have $z^3_v = 0$. Substituting $z^1, z^3$
in $a z = \alpha$, yields $a_{e^j_1}+ a_{h^{j-1}} = 0$. 
Using a similar line of arguments, we obtain:
\begin{equation}\label{ti2}
a_{e^i_1}= a_{e^i_2} = a_{e^{i-1}_3} =a_{e^{i-1}_4} = - a_{h^{i-1}} \qquad \forall i \in [n]\setminus\{1\}.  
\end{equation}
Now consider a tight point $z^4 \in \S_G$ with $z^4_v = 1$ for all $v \in V$, and construct another tight point $z^5 \in \S_G$ with $z^5_v = 1$ for all $v \in V \setminus \{v^j_5\}$ and $z^5_{v^j_5} = 0$ for some $j \in [n] \setminus \{1\}$. Substituting $z^4, z^5$
in $a z = \alpha$, yields $a_{e^j_1}+ a_{V^j} = 0$. 
Using a similar line of arguments, we obtain:
\begin{align}\label{ti3}
& a_{e^i_1}= a_{e^i_2} = a_{e^{i}_3} =a_{e^{i}_4} = - a_{V^i} \qquad \forall i \in [n-1] \setminus \{1\} \nonumber\\
&a_{e^1_3} =a_{e^1_4} = - a_{V^1}\\
&a_{e^n_1}= a_{e^n_2} = -a_{V^n}. \nonumber
\end{align}
Combining~\eqref{ti2} and~\eqref{ti3} and using the fact that $z^4$ defined above is a tight point of inequality~\eqref{densefacet}, we obtain: 
\begin{equation}\label{ti4}
a_{e^i_1} = a_{e^i_2} = a_{e^j_3} = a_{e^j_4} = -a_{h^j} = - a_{V^k} = \frac{\alpha}{2n-3} \qquad \forall i \in [n]\setminus \{1\}, j \in [n-1], k \in [n].  
\end{equation}
Since $a z \leq \alpha$ is nontrivial and valid for $\S_G$, we have $\alpha > 0$. Hence, by~\eqref{ti4}, we conclude that inequality~\eqref{densefacet} coincides with $a z \leq \alpha$ up to a positive scaling implying that it defines a facet of $\MP_G$.
\end{prf}

Notice that the hypergraph $G$ in Proposition~\ref{prop dense} has a fixed rank $r = 8$, while $|E| = 6 n-5$ for all $n \geq 2$.

On the positive side, as a corollary to our main results, we obtain an interesting property of the coefficients in facet-defining inequalities for the multilinear poytope of $\beta$-acyclic hypergraphs.

\begin{corollary}
\label{cor coefficient sum}
Let $G=(V,E)$ be a $\beta$-acyclic hypergraph and let $a z \le b$ be a facet-defining inequality of $\MP_G$ different from $z_p \ge 0$, for $p \in V \cup E$.
Then, $\sum_{p \in V \cup E} a_p = b$.
\end{corollary}

\begin{prf}
Since $G$ is $\beta$-acyclic, by~\cref{th beta iff}, it has a nest point sequence of length $|V|$, say $v_1, \dots, v_n$.
Let $G' = (V, E')$ be the expansion of $G$ w.r.t.~$v_1, \dots, v_n$.
From \cref{lem expansion}, $G'$ is expanded w.r.t.~$v_1, \dots, v_n$.
We then apply \cref{th2} to $G'$. 
Denote by $\bar E$ the set of maximal edges of $G'$, and note that $\bar E \subseteq E$.
System~\eqref{r1}--\eqref{r7} contains nonnegativity constraints for the edges in $\bar E$, and every other inequality $c x \le d$ satisfies $\sum_{p \in V \cup E'} c_p = d$.
The extended variables correspond to the edges in $E' \setminus E$, and none of them are in $\bar E$.
The inequality $a z \le b$ is then obtained from System~\eqref{r1}--\eqref{r7}, by projecting out all the variables in $E' \setminus E$ via Fourier-Motzkin elimination.
The projection consists of nonnegativity constraints on the edges in $\bar E$, and of inequalities that are sums of constraints of the form $c x \le d$ with $\sum_{p \in V \cup E'} c_p = d$.
Since $a z \le b$ is not a nonnegativity constraint, we have $\sum_{p \in V \cup E} a_p = b$.
\end{prf}

We remark that, using Proposition~6 in \cite{dPKha17MOR}, \cref{cor coefficient sum} also holds for facet-defining inequalities of $\MP_G$, for general a hypergraph $G$, provided that their support hypergraphs are $\beta$-acyclic.

%\adp{Move this this corollary earlier? When discussing FM elimination.}

\subsection{Extended running intersection inequalities}
\label{extendedrunning}
Let us consider again the description for the multilinear polytope of an expansion of a $\beta$-acyclic hypergraph $G$ given by inequalities~\eqref{r1}--\eqref{r7}.
Inequalities~\eqref{r1}--\eqref{r4} and inequalities~\eqref{r6}--\eqref{r7} are either flower inequalities or are present in the standard linearization. Now consider inequalities~\eqref{r5}: these inequalities are running intersection inequalities, if $f'(e)$ is a node of $G$, but are not implied by any previously known inequalities for multilinear sets, if $f'(e)$ is an edge of $G$.
Motivated by this observation, we next introduce a new class of cutting planes for multilinear sets that serve as a generalization of running intersection inequalities, introduced in~\cite{dPKha21MOR}. 

In order to define the new inequalities, we first introduce the notion of running intersection property~\cite{BeFaMaYa83}.
A set $F$ of subsets of a finite set $V$ has the \emph{running intersection property} if there exists an ordering $p_1, p_2, \ldots, p_m$ of the sets in $F$ such that
\begin{equation}
\label{ripeq}
\text{for each $k = 2, \dots,m$, there exists $j < k$ such that $p_k \cap \Big(\bigcup_{i < k}{p_i}\Big) \subseteq p_j$.}
\end{equation}
Henceforth, we refer to an ordering $p_1, p_2, \ldots, p_m$ satisfying~\eqref{ripeq} as a \emph{running intersection ordering} of $F$.
Each running intersection ordering $p_1, p_2, \ldots, p_m$ of $F$ induces a collection of sets
\begin{equation}
\label{ripeqsets}
N(p_1) := \emptyset, \qquad N(p_k) := p_k \cap \Big(\bigcup_{i < k}{p_i}\Big) \text{ for } k=2,\dots,m.
\end{equation}

\begin{definition}
\label{def:rie}
Consider a hypergraph $G = (V, E)$.
Let $e_0 \in E$ and let $e_k$, $k \in K$, be a collection of edges in $E$ with $e_0 \cap e_k \neq \emptyset$ for all $k \in K$, such that
the set $\tilde E := \{e_0 \cap e_k : k \in K\}$ has the running intersection property.
%\begin{itemize}
%\item[(i)] $|e_0 \cap e_k| \geq 2$ for all $k \in K$,
%\item[(ii)] $e_0 \cap e_k \not\subseteq e_0 \cap e_{k'}$ for any $k, k' \in K$,
%\item[(iii)] the set $\tilde E := \{e_0 \cap e_k : k \in K\}$
%\end{itemize}
Consider a running intersection ordering of $\tilde E$ with the corresponding sets $N(e_0 \cap e_k)$, for all $k \in K$, as defined in~\eqref{ripeqsets}.
For each $k \in K$, let $w_k \subseteq N(e_0 \cap e_k)$ such that $w_k \in \{\emptyset\} \cup V \cup E$.
We define an~\emph{extended running intersection inequality} centered at $e_0$ with neighbors $e_k$, $k \in K$ as:
\begin{equation}
\label{eq: rie}
- \sum_{k \in K} z_{w_k}
+ \sum_{v \in e_0 \setminus \bigcup_{k \in K} e_k} z_v
+ \sum_{k \in K}{z_{e_k}}
- z_{e_0} \leq \omega-1,
\end{equation}
where we define $z_{\emptyset} = 0$, and
$$\omega = \Big|e_0 \setminus \bigcup_{k \in K} e_k \Big|+\Big|\Big\{k \in K : N(e_0 \cap e_k) = \emptyset \Big\}\Big|.$$
\end{definition}

We do not include the proof of validity for extended running intersection inequalities, as the proof mirrors the proof of validity for running intersection inequalities (see Proposition~1 in~\cite{dPKha21MOR}). In~\cite{dPKha21MOR}, the
authors prove that the system of all running intersection inequalities centered
at $e_0$ with neighbors $e_k$, $k \in K$, is independent of the running intersection
ordering (see Proposition~2 in~\cite{dPKha21MOR}). The same statement holds for extended running intersection inequalities. 

\begin{remark}\label{redERI}
In the special case where the sets $w_k$ for all $k \in K$ with $N(e_0 \cap e_k) \neq \emptyset$ are nodes of $G$, extended running intersection inequalities simplify to running intersection inequalities introduced in~\cite{dPKha21MOR}. In an even more restrictive setting where $w_k = \emptyset$ for all $k \in K$, extended running intersection inequalities simplify to flower inequalities introduced in~\cite{dPKha18SIOPT}.
%It can be checked that an ERI inequality is redundant, \ie it is implied by other ERI inequalities and the standard linearization, if any of the following conditions is not satisfied:
%Moreover if there exist a collection of nested edges $e_i$, $i \in I$ with the same intersection with the center edge $e_0$, it suffices to consider the smallest $e_i$
%in the collection. 
\end{remark}

We now define the~\emph{extended running intersection relaxation} of the multilinear set $\S_G$, denoted by $\MP^{\rm ERI}_G$, as the polytope obtained by adding to the standard linearization, all possible extended running intersection inequalities of $\S_G$.
For a general hypergraph $G$, many of the extended running intersection inequalities are redundant
for $\MP^{\rm ERI}_G$. The following proposition provides sufficient conditions to identify such redundant inequalities.
 
\begin{proposition}\label{prop:red}
   Consider an extended running intersection inequality centered at $e_0$ with neighbors $e_k$, $k \in K$, as defined by~\eqref{eq: rie}. If this inequality defines a facet of $\MP^{\rm ERI}_G$, then it satisfies the following conditions: 
   \begin{itemize}
\item [(i)] for any $k \neq k' \in K$, we have $e_0 \cap e_k \not\subseteq e_0 \cap e_{k'}$,
\item [(ii)] for each $k \in K$, we have $|e_0 \cap e_k| \geq 2$,
\item [(iii)] for any $k \neq k' \in K$, with $w_k , w_{k'} \in N(e_0 \cap e_k) \cap N(e_0 \cap e_{k'})$, we have $w_k = w_{k'}$.
\item  [(iv)] for each $k \in K$, we have $w_k \not\subset p \subseteq N(e_0 \cap e_k)$ for any $p \in E$.
\end{itemize}
\end{proposition}

\begin{prf}
    The proof of redundancy of an extended running intersection inequality not satisfying one of the conditions~(i)--(iii) follows from the proof of Proposition~3 in~\cite{dPKha21MOR} regarding the redundancy of running intersection inequalities.
    Hence it suffices to show that if an extended running intersection inequality does not satisfy condition~(iv), then it is implied by other inequalities in $\MP^{\rm ERI}_G$.

    Consider an extended running intersection inequality centered at $e_0$ with neighbors $e_k$, $k \in K$ such that for some $\bar k \in K$ we have $w_{\bar k} \subset p \subseteq N(e_0 \cap e_{\bar k})$ for some $p \in E$. Then consider another 
    extended running intersection inequality that is identical to the first one except for $w_k$ replaced by $p$. Moreover, consider the inequality $z_p \leq z_{w_{\bar k}}$ present in the standard linearization and hence present in $\MP^{\rm ERI}_G$. Summing up the latter two inequalities, we obtain the first extended running intersection inequality, and this completes the proof.
\end{prf}

Condition~(iv) of Proposition~\ref{prop:red} identifies conditions under which running intersection inequalities are implied by extended running intersection inequalities. The following example demonstrates this fact.

\begin{figure}[h]
\begin{center}
\includegraphics[width=.7\textwidth]{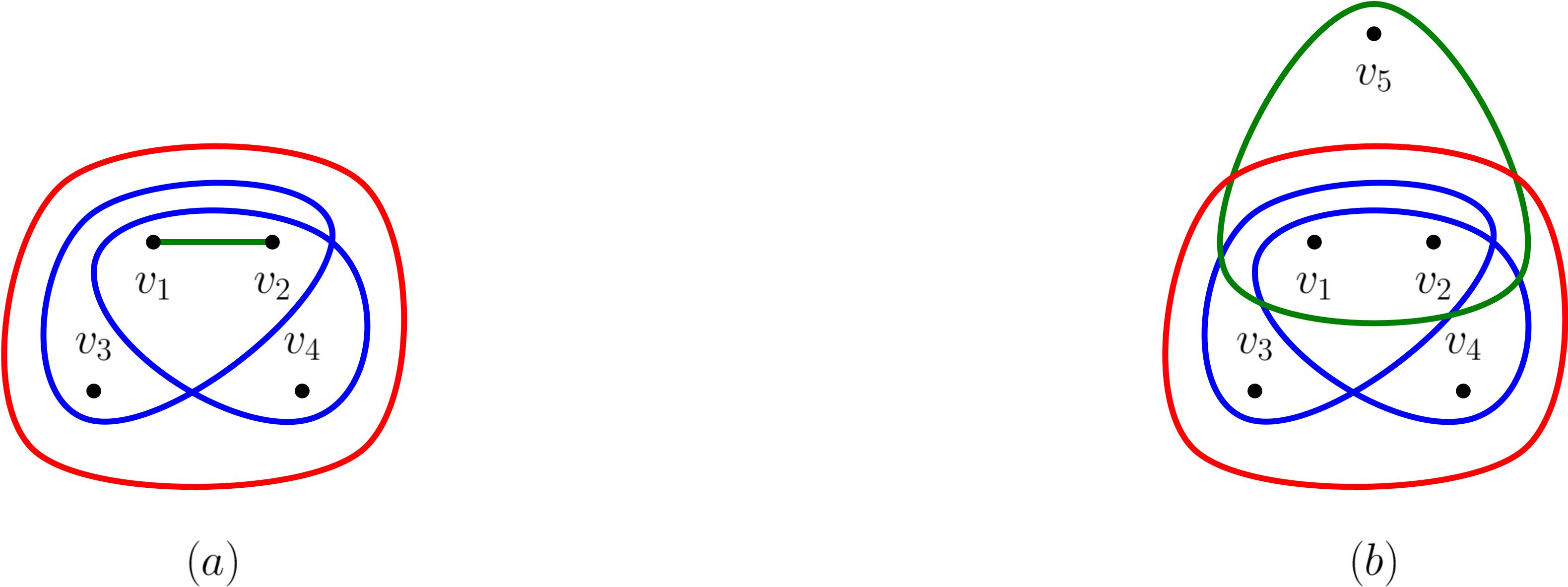}
\end{center}
\caption{Illustration of the hypergraphs considered in \cref{ex hypergraph a,ex hypergraph b}.}
\label{fig twohypergraphs}
\end{figure}

\begin{example}
\label{ex hypergraph a}
Consider the hypergraph $G = (V, E)$ with $V = \{v_1, v_2, v_3, v_4\}$ and 
$$
E=\Big\{\{v_1, v_2\}, \{v_1, v_2, v_3\}, \{v_1, v_2, v_4\}, V\Big\}.
$$ 
See \cref{fig twohypergraphs} $(a)$. 
The running intersection inequalities centered at $V$ with neighbors $\{v_1, v_2, v_3\}$, $\{v_1, v_2, v_4\}$ are given by:
\begin{align}
\label{rIs}
\begin{split}
&-z_1 + z_{123}+z_{124} - z_{1234} \leq 0 \\
& -z_2 + z_{123}+z_{124} - z_{1234} \leq 0.
\end{split}
\end{align}
Moreover, the additional extended running intersection inequality centered at $V$ with neighbors $\{v_1, v_2, v_3\}$, $\{v_1, v_2, v_4\}$ is given by:
\begin{equation}\label{erI}
-z_{12}+ z_{123}+z_{124} - z_{1234} \leq 0.
\end{equation}
Since $z_{12} \leq z_1$ and $z_{12} \leq z_2$, inequality~\eqref{erI} implies inequalities~\eqref{rIs}. It can be checked that inequality~\eqref{erI} defines a facet of $\MP_G$. In fact by adding inequality~\eqref{erI} together with flower inequalities to $\MP^{\rm LP}_{G}$, we obtain $\MP_G$.
\end{example}

 Now suppose that $G$ is $\beta$-acyclic.  Notice that extended running intersection inequalities are sparse; that is, for a rank $r$ hypergraph, extended running intersection inequalities contain at most $2(r-1)$ nonzero coefficients, implying by Proposition~\ref{prop dense} that $\MP^{\rm ERI}_G$ does not coincide with the multilinear polytope of $\beta$-acyclic hypergraphs. We leave as an open question the problem of characterizing the class of hypergraphs $G$ for which we have $\MP_G = \MP^{\rm ERI}_G$.
The following example provides perhaps the simplest $\beta$-acyclic hypergraph $G$ for which we have 
$\MP_G \subset \MP^{\rm ERI}_G$.

\begin{example}
\label{ex hypergraph b}
Consider the $\beta$-acyclic hypergraph $G = (V, E)$ with $V = \{v_1, \cdots, v_5\}$ and 
$$
E=\Big\{\{v_1, v_2, v_3\}, \{v_1, v_2, v_4\}, \{v_1, v_2, v_5\}, 
\{v_1, v_2, v_3, v_4\} \Big\}.
$$ 
See \cref{fig twohypergraphs} $(b)$. 
It can be checked that the following inequality
defines a facet of $\MP_G$:
\begin{equation}\label{notERI}
z_5-z_{125}+z_{123}+z_{124}-z_{1234} \leq 1.
\end{equation}
However, the above inequality is not an extended running intersection inequality since 
$\{v_1, v_2, v_5\} \notin \{v_1, v_2, v_3\} \cap \{v_1, v_2, v_4\}$. In fact, inequality~\eqref{notERI} can be obtained as follows. 
Let $G'=(V, E')$ denote the expansion of $G$ with respect to the nested sequence $\{v_5, v_4, v_3, v_2, v_1\}$. Then we have $E' = E \cup \{\{v_1, v_2\}\}$. By Theorem~\ref{th2}, the following inequalities are implied by $\MP_{G'}$:
\begin{align*}
-z_{12} + z_{123}+z_{124} - z_{1234} \leq 0 \\
z_5 + z_{12} - z_{125} \leq 1.
\end{align*}
Projecting out $z_{12}$ from above inequalities, we obtain inequality~\eqref{notERI}.  Indeed employing this technique in a recursive manner, one can obtain dense facet-defining inequalities for the multilinear polytope of $\beta$-acyclic hypergraphs.
\end{example}

\bigskip
\noindent
\textbf{Acknowledgements:}
The authors would like to thank Silvia Di Gregorio for discussions and preliminary work on the characterization of the multilinear polytope for $\beta$-acyclic hypergraphs.

\bigskip
\noindent
\textbf{Funding:}
A. Del Pia is partially funded by AFOSR grant FA9550-23-1-0433. 
A. Khajavirad is in part supported by AFOSR grant FA9550-23-1-0123.
Any opinions, findings, and conclusions or recommendations expressed in this material are those of the authors and do not necessarily reflect the views of the Air Force Office of Scientific Research.

\ifthenelse {\boolean{MPA}}
{
% For MPA begin
\bibliographystyle{spmpsci}
% For MPA end
}
{
% For OO begin
%\bibliographystyle{plain}
\bibliographystyle{plainurl}
% For OO end
}

%\bibliography{biblio}

\begin{thebibliography}{10}

\bibitem{Bal98}
E.~Balas.
\newblock Disjunctive programming: properties of the convex hull of feasible
  points.
\newblock {\em Discrete Applied Mathematics}, 89(1--3):3--44, 1998.

\bibitem{BeFaMaYa83}
C.~Beeri, R.~Fagin, D.~Maier, and M.~Yannakakis.
\newblock On the desirability of acyclic database schemes.
\newblock {\em Journal of the ACM}, 30:479--513, 1983.

\bibitem{BieMun18}
D.~Bienstock and G.~Munoz.
\newblock Lp formulations for polynomial optimization problems.
\newblock {\em SIAM Journal on Optimization}, 28(2):1121--1150, 2018.

\bibitem{bra14}
J.~Brault-Baron.
\newblock Hypergraph acyclicity revisited.
\newblock {\em ACM Computing Surveys}, 49(3):54:1--54:26, 2016.

\bibitem{BucCraRod16}
C.~Buchheim, Y.~Crama, and E.~Rodr\'iguez-Heck.
\newblock Berge-acyclic multilinear 0--1 optimization problems.
\newblock {\em European Journal of Operational Research}, 273(1):102--107,
  2019.

\bibitem{chenSanOkt20}
R.~Chen, S.~Dash, and O.~G{\"u}nl{\"u}k.
\newblock Cardinality constrained multilinear sets.
\newblock In {\em International Symposium on Combinatorial Optimization}, pages
  54--65. Springer, 2020.

\bibitem{Cor01b}
G.~Cornu\'ejols.
\newblock {\em Combinatorial Optimization: Packing and Covering}, volume~74 of
  {\em CBMS-NSF Regional Conference Series in Applied Mathematics}.
\newblock SIAM, 2001.

\bibitem{yc93}
Y.~Crama.
\newblock {Concave extensions for non-linear $\01$ maximization problems}.
\newblock {\em Mathematical Programming}, 61:53--60, 1993.

\bibitem{CraRod16}
Y.~Crama and E.~Rodr\'iguez-Heck.
\newblock A class of valid inequalities for multilinear $0-1$ optimization
  problems.
\newblock {\em Discrete Optimization}, 25:28--47, 2017.

\bibitem{dPDiG21IJO}
A.~Del~Pia and S.~Di~Gregorio.
\newblock Chv\'atal rank in binary polynomial optimization.
\newblock {\em INFORMS Journal on Optimization}, 3(4):315--349, 2021.

\bibitem{delGre22}
A.~Del~Pia and S.~Di~Gregorio.
\newblock On the complexity of binary polynomial optimization over acyclic
  hypergraphs.
\newblock In {\em Proceedings of the 2022 Annual ACM-SIAM Symposium on Discrete
  Algorithms (SODA)}, pages 2684--2699, 2022.

\bibitem{dPDiG23ALG}
A.~Del~Pia and S.~Di~Gregorio.
\newblock On the complexity of binary polynomial optimization over acyclic
  hypergraphs.
\newblock {\em To appear in Algorithmica}, 2022.

\bibitem{dPKha17MOR}
A.~Del~Pia and A.~Khajavirad.
\newblock A polyhedral study of binary polynomial programs.
\newblock {\em Mathematics of Operations Research}, 42(2):389--410, 2017.

\bibitem{dPKha18SIOPT}
A.~Del~Pia and A.~Khajavirad.
\newblock The multilinear polytope for acyclic hypergraphs.
\newblock {\em SIAM Journal on Optimization}, 28(2):1049--1076, 2018.

\bibitem{dPKha18MPA}
A.~Del~Pia and A.~Khajavirad.
\newblock On decomposability of multilinear sets.
\newblock {\em Mathematical Programming, Series A}, 170(2):387--415, 2018.

\bibitem{dPKha21MOR}
A.~Del~Pia and A.~Khajavirad.
\newblock The running intersection relaxation of the multilinear polytope.
\newblock {\em Mathematics of Operations Research}, 46(3):1008--1037, 2021.

\bibitem{dPKhaSah20MPC}
A.~Del~Pia, A.~Khajavirad, and N.~Sahinidis.
\newblock On the impact of running-intersection inequalities for globally
  solving polynomial optimization problems.
\newblock {\em Mathematical Programming Computation}, 12:165--191, 2020.

\bibitem{dPWal22IPCO}
Alberto Del~Pia and M.~Walter.
\newblock Simple odd $\beta$-cycle inequalities for binary polynomial
  optimization.
\newblock In {\em Proceedings of IPCO 2022}, volume 13265 of {\em Lecture Notes
  in Computer Science}, pages 181--194. Springer, 2022.

\bibitem{Dur12}
D.~Duris.
\newblock Some characterizations of $\gamma$ and $\beta$-acyclicity of
  hypergraphs.
\newblock {\em Information Processing Letters}, 112:617--620, 2012.

\bibitem{fagin83}
Ronald Fagin.
\newblock Degrees of acyclicity for hypergraphs and relational database
  schemes.
\newblock {\em Journal of the ACM (JACM)}, 30(3):514--550, 1983.

\bibitem{HojPfeWal19}
C.~Hojny, M.~Pfetsch, and M.~Walter.
\newblock Integrality of linearizations of polynomials over binary variables
  using additional monomials.
\newblock {\em Preprint, arXiv:1911.06894}, 2019.

\bibitem{Aida22}
A.~Khajavirad.
\newblock On the strength of recursive mccormick relaxations for binary
  polynomial optimization.
\newblock {\em Operations Research Letters}, 51(2):146--152, 2023.

\bibitem{IdaNick18}
A.~Khajavirad and N.~V. Sahinidis.
\newblock {A hybrid LP/NLP paradigm for global optimization relaxations}.
\newblock {\em Mathematical Programming Computation}, 10(3):383--421, May 2018.

\bibitem{kim22}
J.~Kim, J.~P. Richard, and M.~Tawarmalani.
\newblock A reciprocity between tree ensemble optimization and multilinear
  optimization.
\newblock {\em Optimization Online,
  {https://optimization-online.org/2022/03/8828/}}, 2022.

\bibitem{Pad89}
M.~Padberg.
\newblock The {B}oolean quadric polytope: Some characteristics, facets and
  relatives.
\newblock {\em Mathematical Programming}, 45(1--3):139--172, 1989.

\bibitem{rothvoss17}
T.~Rothvoss.
\newblock The matching polytope has exponential extension complexity.
\newblock {\em Journal of the ACM (JACM)}, 64(6):1--19, 2017.

\bibitem{SchBookIP}
Alexander Schrijver.
\newblock {\em Theory of Linear and Integer Programming}.
\newblock Wiley, Chichester, 1986.

\bibitem{XuAdaAks20}
Y.~Xu, W.~Adams, and A.~Gupte.
\newblock Polyhedral analysis of symmetric multilinear polynomials over box
  constraints.
\newblock {\em arXiv preprint arXiv:2012.06394}, 2020.

\end{thebibliography}

\end{document}